\documentclass[10pt,twoside]{article}

\usepackage{amsfonts} %MIORTEGA
\usepackage{units, amsmath, color}
\usepackage{pstricks,pst-node,pst-coil,pst-plot}  

\usepackage{graphicx}

\usepackage{amssymb}
\usepackage{latexsym}

\usepackage{tikz}

\usepackage{hyperref}

\makeatletter

\@addtoreset{equation}{section} \makeatother
\newtheorem{theorem}{Theorem}[section]
\newtheorem{remark}[theorem]{Remark}
\newtheorem{lemma}[theorem]{Lemma}
\newtheorem{proposition}[theorem]{Proposition}
\newtheorem{corollary}[theorem]{Corollary}
\newtheorem{definition}[theorem]{Definition}
\newtheorem{example}[theorem]{Example}
\newcommand{\findemo}{\hfill \qed}

\def\1{\mathbf{1}}
\def\:{\lrcorner}
\def\#{\sharp}

\def\d{\delta}
\def\e{\epsilon}

\def\s{\sigma}

\def\x{\otimes}

\def\qed{\ensuremath{\quad\Box\quad}}

\def\inv#1{\raise.1em\hbox to 0pt{$^{-1}$\hss}_{#1}\;}

\def\V{\noindent}
\def\v{\noindent}

\newcommand{\bean}{\begin{eqnarray*}}
\newcommand{\eean}{\end{eqnarray*}}
\newcommand{\benu}{\begin{enumerate}}
\newcommand{\eenu}{\end{enumerate}}
\newcommand{\eea}{\end{eqnarray}}
\newcommand{\bea}{\begin{eqnarray}}

\newcommand{\ovl}{\overline}

\def \beit{\begin{itemize}}
\def \eeit{\end{itemize}}
\def \bey{\begin{eqnarray*}}
\def \eey{\end{eqnarray*}}

\def \bui#1#2{\mathrel{\mathop{\kern 0pt#1}\limits^{#2}}}
\def \buil#1#2{\mathrel{\mathop{\kern 0pt#1}\limits_{#2}}}

\newtheorem{Lemma}{Lemma}

\newcommand{\be}{\begin{equation}}
\newcommand{\ee}{\end{equation}}

\newcommand{\N}{{\mathbb N}}

\newcommand{\R}{{\mathbb R}}

\newcommand{\C}{{\mathbb C}}

\newcommand{\ben}{\begin{enumerate}}
\newcommand{\een}{\end{enumerate}}
\newcommand{\bit}{\begin{itemize}}
\newcommand{\eit}{\end{itemize}}
\newcommand{\edoc}{\end{document}}

\def \findemo{\hfill$\square$\\}

\def \ovl{\overline}

\newcommand{\bdefi}{\begin{definition}}
\newcommand{\btheo}{\begin{theorem}}
\newcommand{\bprop}{\begin{proposition}}
\newcommand{\brema}{\begin{remark}}
\newcommand{\bcoro}{\begin{corollary}}
\newcommand{\blemm}{\begin{lemma}}
\newcommand{\bexam}{\begin{example}}

\newcommand{\edefi}{\end{definition}}
\newcommand{\etheo}{\end{theorem}}
\newcommand{\eprop}{\end{proposition}}
\newcommand{\erema}{\end{remark}}
\newcommand{\ecoro}{\end{corollary}}
\newcommand{\elemm}{\end{lemma}}
\newcommand{\eexam}{\end{example}}

\title{Global solvability of massless Dirac-Maxwell systems}

\begin{document}

\author{Nicolas Ginoux\footnote{IUT de Metz, D\'epartement Informatique, \^Ile du Saulcy, CS 10628, F-57045 Metz Cedex 01, Email: \texttt{nicolas.ginoux@univ-lorraine.fr}}, Olaf M\"uller\footnote{Fakult\"at f\"ur Mathematik, Universit\"at Regensburg, D-93040 Regensburg, Email: \texttt{olaf.mueller@ur.de}}}

\date{\today}
\maketitle

\begin{abstract}
\V We consider the Cauchy problem of massless Dirac-Maxwell equations on an asymptotically flat background and give a global existence and uniqueness theorem for initial values small in an appropriate weighted Sobolev space. The result can be extended via analogous methods to Dirac-Higgs-Yang-Mills theories.
\end{abstract}

\begin{small}{\it Mathematics Subject Classification} (2010):  35Lxx, 35Qxx, 53A30, 53C50, 53C80
\end{small}

\section{Introduction}

Let $(M^n,g)$ be a globally hyperbolic spin manifold endowed with a trivial $U(1)$-principal bundle $\pi\colon E \rightarrow M $. 
Let $A$ be a connection one-form on $\pi$, or equivalently, a $U(1)$-invariant $i\mathbb{R}$-valued one-form on $E$.
We will assume in the following that $M$ is simply-connected and will regard $A$ as a real-valued one-form on $M$.
We denote the standard spinor bundle of $(M,g)$ by $\sigma\colon \Sigma \rightarrow M$, by $\langle\cdot\,,\cdot\rangle$ the pointwise Hermitian inner product on $\sigma$ and by ``$\,\cdot\,$'' the pointwise Clifford multiplication by vector fields or forms on $\sigma$.
Recall that the Levi-Civita connection $\nabla$ on $TM$ induces a metric covariant derivative on $\sigma$ that we also denote by $\nabla$.
That covariant derivative together with $A$ define a new covariant derivative $\nabla^A$ on $\sigma$ via $ \nabla^A_X (\psi) := \nabla_X \psi + iA(X) \psi $ for any vector field $X$ on $M$.
By definition, the Dirac operator associated to $A$ is the Clifford-trace of $\nabla^A$, that is, for any local orthonormal frame $\left(e_j\right)_{1\leq j\leq n}$ of $TM$, we have $D^A := i \sum_{j=0}^n \e_j e_j \cdot \nabla_{e_j}^A$, where $\e_j=g(e_j,e_j)=\pm1$.
Alternatively, we can write $D^A=D-A\cdot$, where $D$ is the standard Dirac operator of $(M,g)$ and is obtained as the Clifford-trace of $\nabla$.\\

The {\bf Dirac-Maxwell Lagrangian density $\mathcal{L}_{DM}$ for $N$ particles of masses $m_1,\ldots,m_N$ and charges $\mathrm{sgn}(\mu_1)\sqrt{|\mu_1|},\ldots,\mathrm{sgn}(\mu_N)\sqrt{|\mu_N|}$} is defined by 
\[\mathcal{L}_{DM} (\psi \oplus A ) := \frac{1}{4} tr(F^A \wedge F^A) + \sum_{l=1}^N \frac{1}{2}  (\langle D^{\mu_l A} \psi^l , \psi^l  \rangle + \langle \psi^l  , D^{\mu_l A} \psi^l   \rangle )  - \sum_{l=1}^N m_l \langle \psi^l, \psi^l \rangle,\]
where $\psi=(\psi^1,\ldots,\psi^N$ is a section of $\bigoplus_{l=1}^N\s$ and $A$ is a real one-form on $M$.
The critical points of the Lagrangian are exactly the preimages of zero under the operator $ P_{DM}$ given by 
$$ P_{DM} (\psi^1 \oplus ... \oplus \psi^N  \oplus A)  = (D^{\mu_1A} \psi^1 - m_1 \psi^1, \ldots, D^{\mu_NA} \psi^N - m_N \psi^N , d^*dA - J_{\psi} ),$$
 where $J_{\psi} (X)  := \sum_{l=1}^N \mu_l \cdot j_{ll}(X)$ and $j_{kl}(X):= \langle X \cdot \psi^k , \psi^l \rangle $. If $\psi^k$ and $\psi^l$ have equal mass and charge, then it is easy to see that $d^* j_{kl} = 0$, thus in particular $J_\psi$ is divergence-free for $(\psi, A) \in P_{DM}^{-1} (0)$. 
In the sequel, we shall call a pair $(\psi=(\psi^1,\ldots,\psi^N),A)$ as above a solution to the {\bf Dirac-Maxwell equation} if $(\psi,A)\in P_{DM}^{-1} (0)$, that is, if 
\[D^{\mu_l A}\psi^l=m_l\psi^l,\;\;l=1,\ldots,N\qquad\textrm{ and }\qquad d^*dA=J_\psi.\] 
The \emph{massless} Dirac-Maxwell equation is the Dirac-Maxwell equation with $m_1=\ldots=m_N=0$.\\

Let us first shortly review the state of the art on this subject.
Considering the fact that the massless Dirac-Maxwell equation is in dimension $4$ conformally invariant, Christodoulou and Choquet-Bruhat \cite{CBC82} show existence of solutions of Dirac-Yang-Mills-Higgs solutions on four-dimensional Minkowski space with initial values small in weighted Sobolev spaces, the weights being induced by rescaling via the conformal Penrose embedding Minkowski space into the Einstein cylinder. 
One could try to apply their result to Maxwell-Dirac Theory, but, as we are going to explain in the next paragraph, the resulting statement is only nonempty if we extend their setting to a system of finitely many massles particles whose total charge is zero.
Psarelli \cite{mP}, in contrast, treats the question of Dirac-Maxwell equations with or without mass on $\R^{1,3}$ (not in terms of connections modelling potentials, but in terms of curvature tensors modelling field strength\footnote{Recall, however, that the Aharanov-Bohm effect shows that rather than the electromagnetic fields, the potentials play the more fundamental role in electrodynamics}), with results of the form: If $C$ is any compact subset of a Cauchy surface $S$ of $\R^{1,3}$ then there is a number $a$ depending on $C$ such that, if some initial values $I$ with (among others) spinor part supported in $C$ have Sobolev norm smaller than $a$, then there is a global solution with initial values $I$. In the massless case, this result is of course strictly weaker than the weighted Sobolev result.\\
Flato, Simon and Taflin \cite{FlatoSimonTaflin87} were the first to show global existence for \emph{massive} Dirac-Maxwell equations on $\mathbb{R}^{1,3}$ via the construction of explicit approximate solutions and for suitable initial data that are not easy to handle.
For initial data sufficiently small in some weighted Sobolev norm in $\R^{1,3}$, it is Georgiev \cite{G} who established the first global existence result for massless or massive Maxwell-Dirac equations.
The core idea of Georgiev's proof is a gauge in which the potential one-form $A$ satisfies $tA_0 + \sum_{i=1}^3 x^j A_j = 0$ in canonical coordinates of Minkowski space, implying that after the usual transformation to a Maxwell-Klein-Gordon problem the equations satisfy Klainerman's null condition. The entire construction uses canonical coordinates of Minkowski space, and whereas it seems likely that the proof can be generalized to spacetime geometries decaying to Minkowski spacetimes in an appropriate sense, the question of global existence in other spacetime geometries remains completely open.
Let us mention however that, using the complete null structure for Dirac-Maxwell equations from \cite{DAnconaFoschiSelberg2010}, D'Ancona and Selberg can prove \cite{DAnconaSelberg2011} global existence and well-posedness for Dirac-Maxwell equations on $\mathbb{R}^{1,2}$.
The analysis of Dirac-Maxwell equations also includes refining decay estimates, see for instance \cite{BieriMiaoShahshahani14} where the authors show peeling estimates for non-zero-charge Dirac-Klein-Gordon equations with small initial data on $\mathbb{R}^{1,3}$.\\

The aim of the present article is to generalize Georgiev's results to the much more general case of so-called {\em conformally extendible} spacetimes. This latter notion, explained in greater detail in the next section, is located between between asymptotic simplicity and weak asymptotic simplicity and does not require any asymptotics of the curvature tensor along hypersurfaces. Actually, it is very easy to construct examples by hand of conformally extendible manifolds that are not asymptotically flat. Conversely, maximal Cauchy developments of initial values in a weighted Sobolev neighborhood of initial values are known to possess conformal extensions due to criteria developped by  Friedrich and Chrusciel.
%, and we want to show global existence of solutions to the Maxwell-Dirac equation for initial values small in some weighted Sobolev norm.

\bigskip

Our main result is well-posedness of the Cauchy problem for small Lorenz-gauge constrained initial values for massless Dirac-Maxwell systems of vanishing total charge.
A precise formulation is given in the next section. Our method also applies to other field equations, as long as they display an appropriate conformal behaviour and are gauge-equivalent to a semilinear symmetric hyperbolic system admitting a global solution (cf. Appendix). In particular, Dirac-Higgs-Yang-Mills systems as in Choquet-Bruhat's and Christodoulou's article can be handled similarly. The method --- a special sort of ``causal induction'' --- can be found in Section \ref{proofmain} and seems to be completely new.

\bigskip

In a subsequent work, we will furthermore examine the question whether the solutions of the constraint equations of fixed regularity intersected with any open ball around $0$ always form an infinite-dimensional Banach manifold.

\bigskip

\v The article is structured as follows: The second section introduces the concept of conformal extendibility and gives a detailed account of the main result. The third section recalls well-known facts on transformations under which the Dirac-Maxwell equations display some sort of covariance, proves Proposition \ref{obstruction} and derives the constraint equations used in Theorem \ref{main}. The fourth section is devoted to a proof of the main theorem, and the last section is an appendix transferring standard textbook tools for symmetric hyperbolic systems to the case of coefficients of finite (i.e., $C^k$) regularity needed here, a result that should not surprise experts on the fields and for which we do not claim originality by any means.

\bigskip

{\bf Acknowledgements:} It is our pleasure to thank Helmut Abels, Bernd Ammann, Yvonne Choquet-Bruhat, Piotr Chru\'sciel, Felix Finster, Hans Lindblad, Maria Psarelli and Andr\'as Vasy for fruitful discussions and their interest in this work.

\section{The notion of conformal extendibility and the precise statement of the result}

Let us first review some geometric notions as well as introduce some new terminology.\\
A continuous piecewise $C^1$ curve $c$ in a time-oriented Lorentzian manifold $P$ is called {\bf future} if and only if $c'$ is causal future on the $C^1$ pieces, a subset $A$ of $P$ {\bf causally convex} if any causal curve intersects $A$ in the image of a (possibly empty) interval.
A subset $S$ of $P$ is called {\bf Cauchy surface} if and only if any $C^0$-inextendible causal future curve intersects $S$ exactly once, a subset $A$ {\bf future compact} if and only if for any Cauchy surface $S$ of $P$, the subset $J^+ (S) \cap \overline{A} $ is compact.\\

Let $(M,g)$ and $(N,h)$ be globally hyperbolic Lorentzian manifolds, where $g,h$ are supposed to be $C^k$ metrics for some $k\in\mathbb{N}\setminus\{0\}$ (this reduced regularity is essential for our purposes!).
An open conformal embedding $f \in C^k(M,N)$ is said to {\bf $C^k$-extend $g$ conformally} or to be a {\bf $C^k-$conformal extension of $(M,g)$} if and only if $\overline{f(M)}$ is causally convex and future compact.
A globally hyperbolic manifold $(M,g)$ is, called {\bf $C^k$-extendible} for $k \in \N \cup \{ \infty \} $ if and only if there is a $C^k$-conformal extension of $(M,g)$ into a globally hyperbolic manifold.

Whereas Choquet-Bruhat and Christodoulou work with the Penrose embedding which is a $C^{\infty}$-conformal extension of the entire spacetime, it turns out that, in order to generalize the result by Choquet-Bruhat and Christodoulou, we have to generalize our notion of conformal compactification in a twofold way. First, only the timelike future of a Cauchy surface will be conformally embeddable with open image; furthermore, we have to relax the required regularity of the metric of the target manifolds from $C^{\infty}$ to $C^k$. The reason for the second generalization is that we want to include maximal Cauchy developments $(g, \Phi)$ of initial values for Einstein-Klein-Gordon theories that satisfy decay conditions at spatial infinity only for finitely many derivatives (controlled by a single weighted Sobolev norm).
Thus one cannot control higher derivatives at future null infinity. 
Therefore, we need to show a version of the usual existence theorem for symmetric hyperbolic systems for coefficients of finite regularity, which is done in the appendix \ref{s:appendix}.

The second need for modification comes from the fact that the extension via the Penrose embedding into the Einstein cylinder can, of course, be generalized in a straightforward manner to every compact perturbation of the Minkowski metric. But compact perturbations of Minkowski metric are physically rather unrealistic, as (with interactions like Maxwell theory satisfying the dominant energy condition) a nonzero energy-momentum tensor necessarily entails a positive mass of the metric. A positive mass of the metric, in turn, is an obstacle to a smooth extension at spacelike infinity $i_0$, for a discussion see \cite[pp. 180-181]{MasonNicolas}.
Thus we necessarily have a singularity in the surrounding metric at $i_0$, so that we have to restrict to the timelike future of a fixed Cauchy surface.

\bigskip
Results by Anderson and Chru\'sciel (cf. \cite[Theorems 5.2, 6.1 \& 6.2]{AndersonChrusciel}), improving earlier results by Friedrich \cite{Friedrich86} imply that, apart from the --- physically less interesting --- class of compact perturbations of Minkowski space, there is a rich and more realistic class of manifolds which is $C^4$-extendible in the sense above, namely the class of all static initial values with Schwarzschildian ends and small initial values in an appropriate Sobolev space --- see also Corvino's article on this topic \cite{Cor}.
This space of initial values is quite rich, which can be seen by the conformal gluing technique of Corvino and Schoen \cite{CorvinoSchoen}.
This holds in any even dimension. And in the case of a four-dimensional spacetime, there is, in fact, an even larger class of initial values satisfying the conditions of our global existence theorem which is given by a smallness condition to the Einstein initial values in a weighted Sobolev space encoding a good asymptotic decay towards Schwarzschild initial data, cf. the remark following Theorem 6.2 in \cite{AndersonChrusciel} and the remarks following Theorem 2.6 in \cite{Dain07}.
The maximal Cauchy development of any such initial data set carries even a Cauchy temporal function $t$ such that, for all level sets $S_a:=t^{-1}(\{a\})$ of $t$, both $I^\pm (S_a) $ are $C^4$-extendible and thus satisfy even the stronger assumption of Theorem \ref{corollary}.\footnote{This is a remarkable fact as it is a first approach to the question whether {\em Einstein-Dirac-Maxwell theory is stable around zero}, as the stability theorems imply that Einstein-Maxwell theory is stable around zero initial values for 
given small Dirac fields, and our main result implies that Maxwell-Dirac Theory is stable around zero for maximal Cauchy developments of small Einstein 
initial values.}  

%We call a globally hyperbolic manifold $(M,g)$ {\bf synoptic} if and only if for any $p,q \in M$, $I^+(p) \cap I^+(q) \neq \emptyset $. It is easy to see that this is equivalent to the existence of a $C^0$-inextendible smooth timelike future curve $k: [0, \infty ) \rightarrow M$ with $I^- (k([0, \infty)) = M$, cf \cite{O:Hor}. Minkowski spacetime is synoptic whereas Schwarzschild spacetime is not synoptic because of the presence of the horizon. Now, using the structure result for maximal Cauchy developments of Einstein-scalar field theory by Lindblad-Rodnianski \cite{LR} it is easy to see that 

%\btheo
%For any initial value $I$ small in a certain weighted Banach space, the maximal Cauchy development of $I$ is synoptic.  
%\etheo

 %\v {\bf Proof.} It is straightforward to show that the past of the image of the curve $k(s) := (s,0,0,0)$ in cylindrical coordinates $(t,r, \phi, \theta)$ contains every point for the metric $g_{1,3} + h^0 +h^1$ in Lindblad-Rodnianskis terminology, as for any point $p:= (x,y,0,0)$, a past curve from $(x+3y , 0 ,0, 0)$ is defined by $c(s) = (x+3y- ys , ys , 0,0) $ to $p$. The spherical coordinates can be corrected in an arbitrary small time at the beginning of the curve, \hfill \qed

%\bigskip

%\v In the main theorem we will suppose synopticity, but it is a merely technical assumption. As a recipe to obtain a more general result, replace in the proof of the main result the sets the sets $I^-(k(t_i)) $ by $I^-(C_i)$ for a compact exhaustion $\{ C_i \vert i \in \N \}$ of $M$.  

\bigskip

The central insight presented in this article is that the above mentioned weakened notion of conformal extension suffices to establish  --- however slightly less explicit --- weighted Sobolev spaces of initial values allowing for a global solution.
In particular, we do not impose asymptotic flatness: the theorem is, e.g., applicable to any precompact open subset of de Sitter spacetime whose closure is causally convex. 
In order to formulate the main theorem, we need to introduce the constraint equations arising from the transformation of the Dirac-Maxwell equations into a symmetric hyperbolic system.
Since we shall consider \emph{conformal} embeddings of an open subset of the original spacetime $(M,g)$ into another spacetime $(N,h)$, we must fix a Cauchy hypersurface $S$ of $N$ as well as a Cauchy time function $t$ on $N$ with $t^{-1}\left(\{0\}\right)=S$.
Denoting by $h=-\beta dt^2 + g_t$ the induced metric splitting and by $S_\tau:=t^{-1}\left(\{\tau\}\right)$, we let $A_0,A_1\in\Gamma(T^*M_{|_{S_0}})$ and $\psi_0^l\in\Gamma(\sigma_{|_{S_0}})$, $1\leq l\leq N$, be initial data for the Dirac-Maxwell equations.
We call \emph{constraint equations} for $A_0,A_1,\psi_0^l$ the following identities:
\begin{equation}\label{constreq} 
0=\frac{1}{\beta}A_1(\frac{\partial}{\partial t})-\sum_{j=1}^3(\nabla_{e_j}A_0)(e_j)
\end{equation} 
and
\begin{eqnarray}\label{constreq2}
0&=&-(\nabla^{\tan})^*\nabla^{\tan}A_0(\frac{\partial}{\partial t})-\sum_{j=1}^{3}\nabla_{e_j}A_1(e_j)-\frac{1}{2\beta}\mathrm{tr}_{g_t}(\frac{\partial g_t}{\partial t})A_1(\frac{\partial}{\partial t})+\frac{1}{\beta}A_1(\mathrm{grad}_{g_t}(\beta(t,\cdot))) \nonumber \\
& &+\frac{1}{2\beta}\nabla_{\mathrm{grad}_{g_t}(\beta(t,\cdot))}A_0(\frac{\partial}{\partial t})+\frac{1}{2}g_t(\nabla^{\tan}A_0,\frac{\partial g_t}{\partial t})+\mathrm{ric}^M(\frac{\partial}{\partial t},A_0^\sharp)+\sum_{l=1}^N\mu_l j_{\psi_0^l}(\frac{\partial}{\partial t}),
\end{eqnarray} 
where $\left(e_j\right)_j$ is a local $h$-orthonormal basis of $TM$, $(\nabla^{\tan})^*\nabla^{\tan}:=\sum_{j=1}^{n-1}\nabla_{\nabla_{e_j}^{S_t}e_j}-\nabla_{e_j}\nabla_{e_j}$ and the spinors for two conformally related metrics are identified as usual.\\

Every solution in Lorenz gauge, when restricted to a Cauchy hypersurface, satisfies the constraint equation (see Proposition \ref{p:constrLorenz}). Our main theorem is that, conversely, small constrained initial values can be extended to global solutions:
%More explicitly, our main theorem is the following: 

\btheo[Main theorem]
\label{main}
Let $(M,g)$ be a $4$-dimensional globally hyperbolic spacetime with a Cauchy hypersurface $S'$ such that $I^+ (S') $ is $C^4$-extendible in a globally hyperbolic spacetime $(N,h)$.
%and let $u\in C^\infty(M)$ with $e^{2u}g=h$.
Let $P_{DM}$ be the massless Dirac-Maxwell operator for a finite number of fermion fields.
% whose total charge vanishes\footnote{or, equivalently, of a finite number of fermion fields and prescribed charged sources whose total charge vanishes along .}.
Then, for any Cauchy hypersurface $S \subset I^+ (S') $ of $(M,g)$, there is a weighted $W^{4, \infty}$-neighborhood $U$ of $0$ in $\pi \vert_S $  such that  for every initial value $\left(A_0=A_{|_{S_0}},A_1=\frac{\nabla A}{\partial t}_{|_{S_0}},\psi_0^l=\psi_{|_{S_0}}^l\right)$ in $U$ with zero total charge w.r.t. $S$ and satisfying the constraint equations {\rm(\ref{constreq})} and {\rm(\ref{constreq2})} there is a solution $(\psi,A)$ of $P_{DM}(\psi,A)=0$ in all of $I^+(S)$. The weight is explicitly computable from the geometry. 

% \begin{eqnarray}
% \label{constreq} 
% 0&=&\frac{1}{\beta}A_1(\frac{\partial}{\partial t})-\sum_{j=1}^3(\nabla_{e_j}A_0)(e_j)\\
% 0&=&-(\nabla^{\tan})^*\nabla^{\tan}A_0(\frac{\partial}{\partial t})-\sum_{j=1}^{3}\nabla_{e_j}A_1(e_j)-\frac{1}{2\beta}\mathrm{tr}_{g_t}(\frac{\partial g_t}{\partial t})A_1(\frac{\partial}{\partial t})+\frac{1}{\beta}A_1(\mathrm{grad}_{g_t}(\beta(t,\cdot))) \nonumber \\
% & &+\frac{1}{2\beta}\nabla_{\mathrm{grad}_{g_t}(\beta(t,\cdot))}A_0(\frac{\partial}{\partial t})+\frac{1}{2}g_t(\nabla^{\tan}A_0,\frac{\partial g_t}{\partial t})+\mathrm{ric}^M(\frac{\partial}{\partial t},A_0^\sharp)+\sum_{l=1}^N\mu_l j_{\psi_0^l}(\frac{\partial}{\partial t}),
% \end{eqnarray} 
% where $(\nabla^{\tan})^*\nabla^{\tan}:=\sum_{j=1}^{n-1}\nabla_{\nabla_{e_j}^{S_t}e_j}-\nabla_{e_j}\nabla_{e_j}$, $A_0:=A_{|_{S_0}}\in\Gamma(T^*M_{|_{S_0}})$, $A_1:=\frac{\nabla A}{\partial t}_{|_{S_0}}\in\Gamma(T^*M_{|_{S_0}})$ and $\psi_0^l:=\psi_{|_{S_0}}^l\in\Gamma(\sigma_{|_{S_0}}) , $
% where the metric $h= e^{2u}g = -\beta dt^2 + g_t$ is the one of the extension to which also the orthonormal basis $(e_j)_j$ refers, $t$ a time function with $S= S_0 = t^{-1} (\{ 0 \} )$ and the spinors for two conformally related metrics are identified as usual.
\etheo
%{\red\bf Check the notations of the spinor bundle, it seems not to be the same everywhere in the text.}
{\bf Remark 1:} The result and its proof still work if we replace the Dirac-Maxwell system by a general Dirac-Higgs-Yang-Mills systems in the sense of Choquet-Bruhat and Christodoulou, if the Yang-Mills group $G$ is a product of a compact semisimple group and an abelian group and if the Yang-Mills $G$-principal bundle is trivial. 

{\bf Remark 2:} In case $\beta=1$, which can be assumed without loss of generality by the existence of Fermi coordinates w.r.t. $h$ in a neighbourhood of $S$, the constraint equations (\ref{constreq}) and (\ref{constreq2}) simplify to
\begin{eqnarray*} 
0&=&\frac{\partial}{\partial t}\left(A(\frac{\partial}{\partial t})\right)+d^*_S(A_S)+(n-1)H\cdot A(\frac{\partial}{\partial t})\\
0&=&-\Delta_S\left(A(\frac{\partial}{\partial t})\right)+d^*_S\left(\frac{\nabla A}{\partial t}_{S}\right)-3g_t(\nabla^S A_S,W)+A(d^*_S W)+2|W|^2A(\frac{\partial}{\partial t})+\mathrm{ric}^M\left(\frac{\partial}{\partial t},A^\sharp\right)\\
&&+\sum_{l=1}^N\mu_l j_{\psi_0^l}(\frac{\partial}{\partial t}),
\end{eqnarray*}
where $A_S:=\iota_S^*A\in\Gamma(T^*S)$, $\frac{\nabla A}{\partial t}_{S}:=\iota_S^*\frac{\nabla A}{\partial t}\in\Gamma(T^*S)$, $W:=\frac{1}{2}g_t^{-1}\frac{\partial g_t}{\partial t}$ is the Weingarten map of $\iota_S:S\hookrightarrow M$, $H:=\frac{1}{n-1}\mathrm{tr}(W)$ is its mean curvature and $\psi_0^l:=\psi^l_{|_{S}}\in\Gamma(\sigma_{|_{S}})$.\\
%{\red\bf The notations are not compatible with those of the main theorem.}
{\bf Remark 3:} An inspection of the proof shows that the assumption of $C^4$-extendibility of $I^+(S)$ could be replaced by the weaker assumption of weak $C^4$-extendibility, defined as follows: A globally hyperbolic manifold $(A,k)$ is {\bf weakly $C^l$-extendible} if there is a sequence of smooth spacelike hypersurfaces (not necessarily Cauchy) of $(A,k)$ such that $S_n \subset I^+ (S_{n+1}) $, $A= \bigcup_{i \in \N} I^+ (S_n)$ and $I^+ (S_n) $ is $C^l$-extendible, for all $n \in \N$. This generalization could be interesting applied to $ (A,k) = I_M^+ (S)$ for an asymptotically flat spacetime $M$ and hyperboloidal subsets $S_n$.

\bigskip

We can derive as an immediate corollary for the case that $M$ has a Cauchy temporal function $t$ all of whose level sets are ``extendible in both directions''.
Here it is important to note that every conformal extension $\mathcal{I}$ induces a pair of constraint equations $C_{\mathcal{I}}$ as above.
Then we obtain:

\btheo
\label{corollary}
Let $(M,g)$ be a $4$-dimensional globally hyperbolic manifold with a Cauchy temporal function $t$ such that for all level sets $S_a:=t^{-1}(\{a\})$ of $t$, $I^\pm (S_a) $ are both $C^4$-extendible by a conformal extension $\mathcal{I}^\pm(a)$.
Then for every Cauchy surface $S$ such that $t\vert_S$ is bounded, and for any initial values satisfying the neutrality and the constraint equations $C_{\mathcal{I}^-(e)}, C_{\mathcal{I}^+(f)}$ for $ e> {\rm sup } \,t(S) $, $f < {\rm inf}\, t(S)$ and small in the respective Sobolev spaces, there is a global solution on $M$ to the massless Dirac-Maxwell system above extending those initial values. \hfill \qed 
\etheo

%\v {\bf Conventions: Hodge Laplacian versus connection Laplacian, Dirac operator with $i$...}

For the physically interested reader, we make a little more precise what would have to be done to connect our setting to proper QED.
First of all, one should build up the $n$-particle space as the vector space generated by exterior products of classical solutions that are totally antisymmetric under permutations of different spinor fields of equal mass and charge to obtain the usual fermionic commutation relations. 
Expanding in a basis of $\mathrm{Span}(\psi^1,\ldots,\psi^N)$ orthonormal w.r.t. the conserved $L^2$-scalar product $(\psi, \phi ) := \int_S j_{\psi, \phi} (\nu)$ (where $\nu$ is the normal vector field to a Cauchy surface $S$), we see we can w.r.o.g. assume that the spinor fields form a $(\cdot\,, \cdot)$-orthogonal system.
If we have initial values at $S$ in appropriate Sobolev spaces satisfying this condition, so will the restrictions of the solution to any other Cauchy surface due to the divergence-freeness of the $j_{\psi , \phi}$.
%Moreover, for the purpose of the probability interpretation in the particle interpretation, we should assume $(\psi^l, \psi^l ) =c_l$ for all $l $, which, at a first glance, seems to be in conflict with the assumed smallness of initial values. Actually, this is not quite true: We can transfer our main result, which is formulated for initial values small in every component, to the case of many particles of unit $L^2$ norm each and require smallness for the potential only, in the following manner: Consider a many-particle solution of the free Dirac equation $((D - m_1) \psi^1 ,\ldots, (D- m_N) \psi^N) = 0 $ with $J_{\psi} = 0$ and then look for less trivial solutions $(A, \psi^1 + \phi^1,\ldots, \psi^N + \phi^N)$ of the full Dirac-Maxwell equations close to $(0, \psi^1,\ldots, \psi^N)$. The resulting equations are $0= D \phi^l - m \phi^l +  \mu_l A\cdot\phi^l$, $\Box A = 2 J_{\phi , \psi} + J_{\phi}$ which can be treated in the same manner as the original Dirac-Maxwell equations by the procedure explained in the proof of Theorem \ref{main}, where $J_{\phi,\psi}:=\sum_{l=1}^N\mu_lj_{\phi^l,\psi^l}$.
%For simplicity, we will nevertheless restrict to initial values small in each argument in the following.
The neutrality condition $\int_S J_{\psi} (\nu) = 0$ is in the case of an orthonormal system of spinors equivalent to the condition $\sum_{l=1}^N \mu_l =0$. Moreover, in that case, $J_{\psi}$ can be seen as the expectation value of the quantum-mechanical Dirac current operator, cf. \cite[Sec. 3]{Finster09}.
In the end, one would also need to quantize the bosonic potential $A$. Furthermore, one should consider the sum of all $n$-particle spaces to include phenomena like particle creation, particle annihilation, and also possibly the Dirac sea.

\bigskip

\bigskip

\section{Invariances of the Dirac-Maxwell equations}\label{s:invariances}

%\v {\bf Proof of Theorem \ref{obstruction}:}

Let us first recall important well-known invariances of the Dirac-Maxwell equation:

\blemm\label{l:eichinv}
Let $(\psi,A)$ be a solution of the Dirac-Maxwell equations on a spin spacetime $(M^n,g)$.
\benu
\item (Gauge invariance) For any $f\in C^\infty(M,\R)$, the pair $(\psi':=(e^{-i\mu_1 f}\psi^1,\ldots,e^{-i\mu_N f}\psi^N),A':=A+df)$ solves again the Dirac-Maxwell equations on $(M^n,g)$.
\item (Conformal invariance) If $n=4$, then for any $u\in C^\infty(M,\R)$, the pair $(\ovl{\varphi}:=e^{-\frac{3}{2}u}\ovl{\psi},A)$ solves $D_{\ovl{g}}^{\mu_l A}\ovl{\varphi^l}=m_le^{-u}\ovl{\varphi^l}$ and $d^*_{\ovl{g}}dA=\sum_{l=1}^n\mu_l j_{\ovl{\varphi^l}}$ on $(M^n,\ovl{g}:=e^{2u}g)$, where $\psi\mapsto\ovl{\psi}$, $S_gM\x E\to S_{\ovl{g}}M\x E$, denotes the natural unitary isomorphism induced by the conformal change of metric.
In particular, in dimension $4$, the Dirac-Maxwell equations are scaling-invariant and the massless Dirac-Maxwell equations are even conformally invariant.
\eenu
\elemm 

\begin{proof}
Both statements follow from elementary computations.
For the sake of simplicity, we perform the proof only for $N=1$ and $q=1$.\\
%{\red\bf CHECK!!}\\
1. By definition of the Dirac operator, we have $D^{A'}=D^A-df\cdot$, 
\begin{eqnarray*} 
D^{A'}\psi'&=&(D^A-df\cdot)(e^{-if}\psi)\\
&=&i\cdot(-ie^{-if}df)\cdot\psi+e^{-if}D^A\psi-e^{-if}df\cdot\psi\\
&=&m\psi'
\end{eqnarray*}  
and $d^*dA'=d^* dA+d^* d^2f=d^* dA=j_\psi=j_{\psi'}$.\\
2. First, we compute, for all tangential vector fields $X,Y,Z$ and every $2$-form $\omega$ on $M^n$:
\begin{eqnarray*} 
(\nabla_X^{\ovl{g}}\omega)(Y,Z)&=&X(\omega(Y,Z))-\omega(\nabla_X^{\ovl{g}}Y,Z)-\omega(Y,\nabla_X^{\ovl{g}}Z)\\
&=&X(\omega(Y,Z))-\omega\left(\nabla_X^gY+X(u)Y+Y(u)X-g(X,Y)\mathrm{grad}_g(u),Z\right)\\
& &-\omega\left(Y,\nabla_X^gZ+X(u)Z+Z(u)X-g(X,Z)\mathrm{grad}_g(u)\right)\\
&=&(\nabla_X^g\omega)(Y,Z)-2X(u)\omega(Y,Z)-Y(u)\omega(X,Z)+Z(u)\omega(X,Y)\\
& &+g(X,Y)\omega(\mathrm{grad}_g(u),Z)-g(X,Z)\omega(\mathrm{grad}_g(u),Y).
\end{eqnarray*} 
We deduce that, for the divergence, we have, in a local $g$-ONB $(e_j)_{0\leq j\leq n-1}$ of $TM$ and for every $X\in \Gamma (M,TM)  $,
\begin{eqnarray*} 
(d^*_{\ovl{g}}\omega)(X)&=&-\sum_{j=0}^{n-1}\varepsilon_j(\nabla_{\ovl{e_j}}^{\ovl{g}}\omega)(\ovl{e_j},X)\\
&=&-e^{-2u}\sum_{j=0}^{n-1}\varepsilon_j(\nabla_{e_j}^{\ovl{g}}\omega)(e_j,X)\\
&=&-e^{-2u}\sum_{j=0}^{n-1}\varepsilon_j\Big((\nabla_{e_j}^g\omega)(e_j,X)-2e_j(u)\omega(e_j,X)-e_j(u)\omega(e_j,X)+X(u)\underbrace{\omega(e_j,e_j)}_{0}\\
& &\phantom{-e^{-2u}\sum_{j=0}^{n-1}\varepsilon_j\Big(}+g(e_j,e_j)\omega(\mathrm{grad}_g(u),X)-g(e_j,X)\omega(\mathrm{grad}_g(u),e_j)\Big)\\
&=&e^{-2u}\Big((d^*_g \omega)(X)-(n-4)\omega(\mathrm{grad}_g(u),X)\Big),
\end{eqnarray*} 
that is, $d^*_{\ovl{g}}\omega=e^{-2u}\big(d^*_g \omega-(n-4)\mathrm{grad}_g(u)\lrcorner\omega\big)$.
If in particular $n=4$, then $d^*_{\ovl{g}}\omega=e^{-2u}d^*_g \omega$, so that $d^*_{\ovl{g}}dA=e^{-2u}d^*_g dA$.
On the other hand, the operator $D^A$ is conformally covariant, that is, $D_{\ovl{g}}^A(e^{-\frac{n-1}{2}u}\ovl{\psi})=e^{-\frac{n+1}{2}u}\ovl{D_g^A\psi}$, in particular we have 
\begin{eqnarray*} 
D_{\ovl{g}}^A\ovl{\varphi}&=&D_{\ovl{g}}^A(e^{-\frac{n-1}{2}u}\ovl{\psi})\\
&=&e^{-\frac{n+1}{2}u}\ovl{D_g^A\psi}\\
&=&-me^{-u}\ovl{\varphi}.
\end{eqnarray*} 
It remains to notice that, for every $X\in TM$,
\begin{eqnarray*} 
j_{\ovl{\varphi}}(X)&=&\langle X\cdot_{\ovl{g}}\ovl{\varphi},\ovl{\varphi}\rangle\\
&=&e^{-(n-1)u}\langle X\cdot_{\ovl{g}}\ovl{\psi},\ovl{\psi}\rangle\\
&=&e^{-(n-1)u}e^u\langle \ovl{X\cdot_{g}\psi},\ovl{\psi}\rangle\\
&=&e^{-(n-2)u}\langle X\cdot_g\psi,\psi\rangle\\
&=&e^{-(n-2)u}j_\psi(X),
\end{eqnarray*} 
that is, $j_{\ovl{\varphi}}=e^{-(n-2)u}j_\psi$.
We deduce that, for $n=4$, we have $d^*_{\ovl{g}}dA=e^{-2u}j_\psi=j_{\ovl{\varphi}}$, which concludes the proof.
\findemo
\end{proof}

%{\bf Define Dirac-wave equation}
The {\bf Dirac-wave operator} $P_{DW}$ is defined by 

$$ P_{DW} (\psi^1 \oplus ... \oplus \psi^N  \oplus A) := (D^A \psi^1 - m_1 \psi^1, \ldots, D^A \psi^N - m_N \psi^N , \Box A - J_{\psi} ),$$
 and the {\bf Dirac-wave equation} is just the equation $P_{DW}(\psi,A)=0$, where $\Box:=dd^*+d^* d$.

\bprop[Lorenz gauge]\label{p:constrLorenz}
Let $(M,g)$ be as above.
\beit\item[i)] For any solution $(\psi,A)$ of the Dirac-wave equation, $\Box(d^* A)=0$ holds on $M$.
In particular $d^* A=0$ on $M$ if and only if $(d^* A)_{|_{S_0}}=0=\left(\frac{\partial}{\partial t}d^* A\right)_{|_{S_0}}$.
\item[ii)] Given any solution $(\psi,A)$ to the Dirac-wave equation, the equations $(d^* A)_{|_{S_0}}=0=\left(\frac{\partial}{\partial t}d^* A\right)_{|_{S_0}}$ are equivalent to

\begin{eqnarray}
%\label{constreq} 
0&=&\frac{1}{\beta}A_1(\frac{\partial}{\partial t})-\sum_{j=1}^3(\nabla_{e_j}A_0)(e_j)\\
0&=&-(\nabla^{\tan})^*\nabla^{\tan}A_0(\frac{\partial}{\partial t})-\sum_{j=1}^{3}\nabla_{e_j}A_1(e_j)-\frac{1}{2\beta}\mathrm{tr}_{g_t}(\frac{\partial g_t}{\partial t})A_1(\frac{\partial}{\partial t})+\frac{1}{\beta}A_1(\mathrm{grad}_{g_t}(\beta(t,\cdot))) \nonumber \\
& &+\frac{1}{2\beta}\nabla_{\mathrm{grad}_{g_t}(\beta(t,\cdot))}A_0(\frac{\partial}{\partial t})+\frac{1}{2}g_t(\nabla^{\tan}A_0,\frac{\partial g_t}{\partial t})+\mathrm{ric}^M(\frac{\partial}{\partial t},A_0^\sharp)+\sum_{l=1}^N \mu_l j_{\psi_0^l}(\frac{\partial}{\partial t}),
\end{eqnarray} 
where $A_0:=A_{|_{S_0}}\in\Gamma(T^*M_{|_{S_0}})$, $A_1:=\frac{\nabla A}{\partial t}_{|_{S_0}}\in\Gamma(T^*M_{|_{S_0}})$ and $\psi_0^l:=\psi^l_{|_{S_0}}\in\Gamma(\sigma_{|_{S_0}})$.
\eeit
\eprop

\begin{proof}
Let $(\psi,A)$ solve the Dirac-wave equation.
Then $\Box(d^* A)=d^*(\Box A)=d^* J_\psi$.
But a direct calculation leads to 
\[d^* j_\psi^k=i\left(\langle D^A\psi^k,\psi^k\rangle-\langle\psi^k,D^A\psi^k\rangle\right)=-2\mathrm{Im}(\langle D^A\psi^k,\psi^k\rangle),\]
hence $d^* J_\psi=0$ as soon as $D^A\psi^k=m_k\psi^k$ with $m_k\in \R$ (or, more generally, if $D^A\psi=H\psi$ for some Hermitian endomorphism-field $H$ of $\sigma$).
This shows $\Box(d^* A)=0$ and $i)$.\\

Next we express the equations $(d^* A)_{|_{S_0}}=0=\left(\frac{\partial}{\partial t}d^* A\right)_{|_{S_0}}$ solely in terms of the initial data $A_0$, $A_1$ and $\psi_0$.
It is already obvious that the first equation $(d^* A)_{|_{S_0}}=0$ only depends on $A_0$ (and its tangential derivatives along $S_0$) and $A_1$, however the second equation $\left(\frac{\partial}{\partial t}d^* A\right)_{|_{S_0}}=0$, which contains a derivative of second order in $t$ of $A$, requires the wave equation $\Box A=J_\psi$ in order to yield a relationship between the initial data.\\
Denoting by $(e_j)_{1\leq j\leq 3}$ a local o.n.b. of $TS_0$ and letting $e_0:=\frac{1}{\sqrt{\beta}}\frac{\partial}{\partial t}$ (the future-oriented unit normal field on $S_0$), we have 
\begin{eqnarray*} 
d^* A&=&-\sum_{j=0}^3\varepsilon_j(\nabla_{e_j}A)(e_j)\\
&=&(\nabla_{e_0}A)(e_0)-\sum_{j=1}^3(\nabla_{e_j}A)(e_j)\\
&=&\frac{1}{\beta}\frac{\nabla A}{\partial t}(\frac{\partial}{\partial t})-\sum_{j=1}^3(\nabla_{e_j}A)(e_j).
\end{eqnarray*}
As a first consequence, if we restrict that identity to $S_0$, we obtain
\[(d^* A)_{|_{S_0}}=\frac{1}{\beta}A_1(\frac{\partial}{\partial t})-\sum_{j=1}^3(\nabla_{e_j}A_0)(e_j).\]
Note here that the second term is in general not the divergence of the pull-back of $A_0$ on $S_0$ since the second fundamental form of $S_0$ in $M$ may be non-vanishing.
Differentiating further, we also obtain
\begin{eqnarray*} 
\frac{\partial}{\partial t}d^* A&=&\frac{\partial}{\partial t}\left(\frac{1}{\beta}\frac{\nabla A}{\partial t}(\frac{\partial}{\partial t})\right)-\sum_{j=1}^3\frac{\partial}{\partial t}\left((\nabla_{e_j}A)(e_j)\right)\\
&=&\frac{1}{\beta}\left\{-\frac{1}{\beta}\frac{\partial\beta}{\partial t}\frac{\nabla A}{\partial t}(\frac{\partial}{\partial t})+\frac{\nabla^2 A}{\partial t^2}(\frac{\partial}{\partial t})+\frac{\nabla A}{\partial t}(\frac{\nabla}{\partial t}\frac{\partial}{\partial t})\right\}\\
& &-\sum_{j=1}^{3}\frac{\nabla}{\partial t}\nabla_{e_j}A(e_j)-\sum_{j=1}^{3}\nabla_{e_j}A(\frac{\nabla e_j}{\partial t}),
\end{eqnarray*}
where 
\begin{eqnarray*} 
\sum_{j=1}^{3}\frac{\nabla}{\partial t}\nabla_{e_j}A(e_j)&=&\sum_{j=1}^{3}\nabla_{e_j}\frac{\nabla A}{\partial t}(e_j)+\nabla_{[\frac{\partial}{\partial t},e_j]}A(e_j)+(R_{\frac{\partial}{\partial t},e_j}A)(e_j)\\
&=&\sum_{j=1}^{3}\nabla_{e_j}\frac{\nabla A}{\partial t}(e_j)+\nabla_{[\frac{\partial}{\partial t},e_j]}A(e_j)-A(R_{\frac{\partial}{\partial t},e_j}e_j)\\
&=&\sum_{j=1}^{3}\nabla_{e_j}\frac{\nabla A}{\partial t}(e_j)+\nabla_{[\frac{\partial}{\partial t},e_j]}A(e_j)-\mathrm{ric}^M(\frac{\partial}{\partial t},A^\sharp).
\end{eqnarray*}
Using the equation $\Box A=J_\psi$, we express $\frac{\nabla^2 A}{\partial t^2}$ in terms of $\psi$ and of tangential (up to second order) and normal (up to first order) derivatives of $A$.
Since the metric $g$ has the form $g=-\beta dt^2\oplus g_t$, we can split the rough d'Alembert operator $\Box^\nabla$ (associated to an arbitrary connection $\nabla$ on the bundle under consideration) under the form
\begin{eqnarray}\label{eq:Boxnabla}
\nonumber\Box^\nabla&=&\sum_{j=0}^{n-1}\varepsilon_j(\nabla_{\nabla_{e_j}^Me_j}-\nabla_{e_j}\nabla_{e_j})\\
\nonumber&=&(\frac{1}{\sqrt{\beta}}\frac{\nabla}{\partial t})^2-\frac{1}{\sqrt{\beta}}\nabla_{\nabla_{\frac{\partial}{\partial t}}^M\frac{1}{\sqrt{\beta}}\frac{\partial}{\partial t}}+\sum_{j=1}^{n-1}\nabla_{\nabla_{e_j}^\perp e_j}+\underbrace{\sum_{j=1}^{n-1}\nabla_{\nabla_{e_j}^{S_t}e_j}-\nabla_{e_j}\nabla_{e_j}}_{=:(\nabla^{\tan})^*\nabla^{\tan}}\\
\nonumber&=&\frac{1}{\beta}(\frac{\nabla}{\partial t})^2+\frac{1}{\sqrt{\beta}}\frac{\partial}{\partial t}(\frac{1}{\sqrt{\beta}})\frac{\nabla}{\partial t}-\frac{1}{\sqrt{\beta}}\frac{\partial}{\partial t}(\frac{1}{\sqrt{\beta}})\frac{\nabla}{\partial t}-\frac{1}{\beta}\nabla_{\nabla_{\frac{\partial}{\partial t}}^M\frac{\partial}{\partial t}}\\
\nonumber& &+\frac{1}{2\beta}\mathrm{tr}_{g_t}(\frac{\partial g_t}{\partial t})\frac{\nabla}{\partial t}+(\nabla^{\tan})^*\nabla^{\tan}\\
\nonumber&=&\frac{1}{\beta}(\frac{\nabla}{\partial t})^2-\frac{1}{2\beta^2}\frac{\partial\beta}{\partial t}\frac{\nabla}{\partial t}-\frac{1}{2\beta}\nabla_{\mathrm{grad}_{g_t}(\beta(t,\cdot))}+\frac{1}{2\beta}\mathrm{tr}_{g_t}(\frac{\partial g_t}{\partial t})\frac{\nabla}{\partial t}+(\nabla^{\tan})^*\nabla^{\tan}\\
&=&\frac{1}{\beta}\Big((\frac{\nabla}{\partial t})^2+\frac{1}{2}\{\mathrm{tr}_{g_t}(\frac{\partial g_t}{\partial t})-\frac{1}{\beta}\frac{\partial\beta}{\partial t}\}\frac{\nabla}{\partial t}\Big)+(\nabla^{\tan})^*\nabla^{\tan}-\frac{1}{2\beta}\nabla_{\mathrm{grad}_{g_t}(\beta(t,\cdot))},
\end{eqnarray} 
where, as usual, $(e_j)_{0\leq j\leq n-1}$ denotes a local ONB of $TM$ with $e_0=\frac{1}{\sqrt{\beta}}\frac{\partial}{\partial t}$ and $\varepsilon_j=g(e_j,e_j)\in\{\pm1\}$, the Levi-Civita connections of $(M,g)$ and $(S,g_t)$ are denoted respectively by $\nabla^M$ and $\nabla^{S_t}$ and where we have made use of the following identities (which are easy to check using Koszul's identity):

\[\nabla_X^MY=\nabla_X^{S_t}Y+\nabla_X^\perp Y=\nabla_X^{S_t}Y+\frac{1}{2\beta}\frac{\partial g_t}{\partial t}(X,Y)\frac{\partial}{\partial t}\]
for all $X,Y\in TS_t=T(\{t\}\times S)$ and 
\[\nabla_{\frac{\partial}{\partial t}}^M\frac{\partial}{\partial t}=\frac{1}{2\beta}\frac{\partial\beta}{\partial t}\frac{\partial}{\partial t}+\frac{1}{2}\mathrm{grad}_{g_t}(\beta(t,\cdot)).\]
As a consequence, (\ref{eq:Boxnabla}) gives
\[\frac{1}{\beta}\frac{\nabla^2 A}{\partial t^2}=\left(\Box-(\nabla^{\tan})^*\nabla^{\tan}\right)A+\frac{1}{2\beta}\nabla_{\mathrm{grad}_{g_t}(\beta(t,\cdot))}A+\frac{1}{2\beta}\left(\frac{1}{\beta}\frac{\partial\beta}{\partial t}-\mathrm{tr}_{g_t}(\frac{\partial g_t}{\partial t})\right)\frac{\nabla A}{\partial t}.\]
If $\Box A=J_\psi$, then we deduce that 
\[\frac{1}{\beta}\frac{\nabla^2 A}{\partial t^2}(\frac{\partial}{\partial t})=J_{\psi}(\frac{\partial}{\partial t})-(\nabla^{\tan})^*\nabla^{\tan}A(\frac{\partial}{\partial t})+\frac{1}{2\beta}\nabla_{\mathrm{grad}_{g_t}(\beta(t,\cdot))}A(\frac{\partial}{\partial t})+\frac{1}{2\beta}\left(\frac{1}{\beta}\frac{\partial\beta}{\partial t}-\mathrm{tr}_{g_t}(\frac{\partial g_t}{\partial t})\right)\frac{\nabla A}{\partial t}(\frac{\partial}{\partial t}).\]
Using again the above identities connecting the Levi-Civita connections of $S_t$ and $M$, we obtain
\begin{eqnarray*} 
\frac{\partial}{\partial t}d^* A&=&-\frac{1}{\beta^2}\frac{\partial\beta}{\partial t}\frac{\nabla A}{\partial t}(\frac{\partial}{\partial t})+\frac{1}{\beta}\frac{\nabla A}{\partial t}(\frac{1}{2\beta}\frac{\partial\beta}{\partial t}\frac{\partial}{\partial t}+\frac{1}{2}\mathrm{grad}_{g_t}(\beta(t,\cdot)))\\
& &+J_\psi(\frac{\partial}{\partial t})-(\nabla^{\tan})^*\nabla^{\tan}A(\frac{\partial}{\partial t})+\frac{1}{2\beta}\nabla_{\mathrm{grad}_{g_t}(\beta(t,\cdot))}A(\frac{\partial}{\partial t})+\frac{1}{2\beta}\left(\frac{1}{\beta}\frac{\partial\beta}{\partial t}-\mathrm{tr}_{g_t}(\frac{\partial g_t}{\partial t})\right)\frac{\nabla A}{\partial t}(\frac{\partial}{\partial t})\\
& &+\mathrm{ric}^M(\frac{\partial}{\partial t},A^\sharp)-\sum_{j=1}^{3}\nabla_{e_j}\frac{\nabla A}{\partial t}(e_j)+\nabla_{[\frac{\partial}{\partial t},e_j]}A(e_j)+\nabla_{e_j}A(\frac{\nabla e_j}{\partial t})\\
&=&-(\nabla^{\tan})^*\nabla^{\tan}A(\frac{\partial}{\partial t})-\sum_{j=1}^{3}\nabla_{e_j}\frac{\nabla A}{\partial t}(e_j)-\frac{1}{2\beta}\mathrm{tr}_{g_t}(\frac{\partial g_t}{\partial t})\frac{\nabla A}{\partial t}(\frac{\partial}{\partial t})+\frac{1}{2\beta}\frac{\nabla A}{\partial t}(\mathrm{grad}_{g_t}(\beta(t,\cdot)))\\
& &+\frac{1}{2\beta}\nabla_{\mathrm{grad}_{g_t}(\beta(t,\cdot))}A(\frac{\partial}{\partial t})-\sum_{j=1}^3\nabla_{[\frac{\partial}{\partial t},e_j]}A(e_j)+\nabla_{e_j}A(\frac{\nabla e_j}{\partial t})+\mathrm{ric}^M(\frac{\partial}{\partial t},A^\sharp)+J_\psi(\frac{\partial}{\partial t}).
\end{eqnarray*} 
Now using $\frac{\nabla g}{\partial t}=0$ as well as $\nabla_{e_j}\frac{\partial}{\partial t}=\frac{1}{2\beta}e_j(\beta)\frac{\partial}{\partial t}+\frac{1}{2}g_t^{-1}\frac{\partial g_t}{\partial t}(e_j,\cdot)$, we have 
\begin{eqnarray*} 
\sum_{j=1}^3\nabla_{[\frac{\partial}{\partial t},e_j]}A(e_j)+\nabla_{e_j}A(\frac{\nabla e_j}{\partial t})&=&\sum_{j=1}^3\nabla_{\frac{\nabla e_j}{\partial t}-\nabla_{e_j}\frac{\partial}{\partial t}}A(e_j)+\nabla_{e_j}A(\frac{\nabla e_j}{\partial t})\\
&=&\underbrace{\sum_{j=1}^3\nabla_{\frac{\nabla e_j}{\partial t}}A(e_j)+\nabla_{e_j}A(\frac{\nabla e_j}{\partial t})}_{0}-\sum_{j=1}^3\nabla_{\nabla_{e_j}\frac{\partial}{\partial t}}A(e_j)\\
&=&-\sum_{j=1}^3\frac{1}{2\beta}e_j(\beta)\frac{\nabla A}{\partial t}(e_j)+\frac{1}{2}\nabla_{g_t^{-1}\frac{\partial g_t}{\partial t}(e_j,\cdot)}A(e_j)\\
&=&-\frac{1}{2\beta}\frac{\nabla A}{\partial t}(\mathrm{grad}_{g_t}(\beta(t,\cdot)))-\frac{1}{2}\sum_{j=1}^3\nabla_{g_t^{-1}\frac{\partial g_t}{\partial t}(e_j,\cdot)}A(e_j)\\
&=&-\frac{1}{2\beta}\frac{\nabla A}{\partial t}(\mathrm{grad}_{g_t}(\beta(t,\cdot)))-\frac{1}{2}g_t(\nabla^{\tan}A,\frac{\partial g_t}{\partial t}),
\end{eqnarray*} 
so that we get
\begin{eqnarray*} 
\frac{\partial}{\partial t}d^* A&=&-(\nabla^{\tan})^*\nabla^{\tan}A(\frac{\partial}{\partial t})-\sum_{j=1}^{3}\nabla_{e_j}\frac{\nabla A}{\partial t}(e_j)-\frac{1}{2\beta}\mathrm{tr}_{g_t}(\frac{\partial g_t}{\partial t})\frac{\nabla A}{\partial t}(\frac{\partial}{\partial t})+\frac{1}{\beta}\frac{\nabla A}{\partial t}(\mathrm{grad}_{g_t}(\beta(t,\cdot)))\\
& &+\frac{1}{2\beta}\nabla_{\mathrm{grad}_{g_t}(\beta(t,\cdot))}A(\frac{\partial}{\partial t})+\frac{1}{2}g_t(\nabla^{\tan}A,\frac{\partial g_t}{\partial t})+\mathrm{ric}^M(\frac{\partial}{\partial t},A^\sharp)+J_\psi(\frac{\partial}{\partial t}).
\end{eqnarray*}
Restricting that equation onto $S_0$, we come to 
\begin{eqnarray*} 
\left(\frac{\partial}{\partial t}d^* A\right)_{|_{S_0}}&=&-(\nabla^{\tan})^*\nabla^{\tan}A_0(\frac{\partial}{\partial t})-\sum_{j=1}^{3}\nabla_{e_j}A_1(e_j)-\frac{1}{2\beta}\mathrm{tr}_{g_t}(\frac{\partial g_t}{\partial t})A_1(\frac{\partial}{\partial t})+\frac{1}{\beta}A_1(\mathrm{grad}_{g_t}(\beta(t,\cdot)))\\
& &+\frac{1}{2\beta}\nabla_{\mathrm{grad}_{g_t}(\beta(t,\cdot))}A_0(\frac{\partial}{\partial t})+\frac{1}{2}g_t(\nabla^{\tan}A_0,\frac{\partial g_t}{\partial t})+\mathrm{ric}^M(\frac{\partial}{\partial t},A_0^\sharp)+J_{\psi_0}(\frac{\partial}{\partial t}).
\end{eqnarray*}
This yields the second equation and concludes the proof.
\hfill\qed\\
\end{proof}

\bprop
\label{obstruction}
Let $(\psi = (\psi^1,\ldots, \psi^N),A)$ be any classical solution to the Dirac-Maxwell equation such that, along a given (smooth, spacelike) Cauchy hypersurface $S$ with future-directed unit normal $\nu$, the $1$-form $dA(\nu,  \cdot)$ is compactly supported.
Then $\int_{S'} J_\psi(\nu') =0$ for all Cauchy hypersurfaces $S'$ of $M$ with future unit normal vector $\nu'$.
In particular, for $N= 1$ and $\mu_1 \neq 0$, we can conclude $\psi_1 = 0$. 
\eprop

\begin{proof}
Let $(\psi,A)$ be any classical (i.e., sufficiently smooth) $1$-particle solution to the Dirac-Maxwell equation, that is, $D^A\psi=m\psi$ and $d^* dA=j_\psi$.
Let $S\subset M$ be any smooth spacelike Cauchy hypersurface and $\nu$ be the future-directed unit normal vector field along $S$.
We first compute the codifferential along $S$ of the $1$-form $\nu\lrcorner dA=dA(\nu,\cdot)$.
Let $\{e_j\}_{1\leq j\leq n-1}$ be any local $g$-orthonormal frame on $S$, then
\begin{eqnarray*} 
d^*_S(\nu\lrcorner dA)&=&-\sum_{j=1}^{n-1}e_j\lrcorner\nabla_{e_j}^S(\nu\lrcorner dA)\\
&=&-\sum_{j=1}^{n-1}e_j\lrcorner\left(\nabla_{e_j}^M(\nu\lrcorner dA)-dA(\nu,\nabla_X^M\nu)\nu^\flat\right)\\
&=&-\sum_{j=1}^{n-1}e_j\lrcorner\nabla_{e_j}^M(\nu\lrcorner dA)\\
&=&-\sum_{j=1}^{n-1}e_j\lrcorner\left((\nabla_{e_j}^MdA)(\nu,\cdot)+dA(\nabla_{e_j}^M\nu,\cdot)\right)\\
&=&-\sum_{j=1}^{n-1}(\nabla_{e_j}^MdA)(\nu,e_j)-\sum_{j=1}^{n-1}dA(\nabla_{e_j}^M\nu,e_j),
\end{eqnarray*} 
where the last sum vanishes since $(X,Y)\mapsto g(\nabla_X^M\nu,Y)$ is symmetric.
We are left with 
\[d^*_S(\nu\lrcorner dA)=-(d^*_M dA)(\nu)=-j_\psi(\nu).\]
As a consequence, if $\nu\lrcorner dA$ has compact support on $S$, then by the divergence theorem,
\[\int_S j_\psi(\nu) d\sigma_g=-\int_S d^*_S(\nu\lrcorner dA) d\sigma_g=0.\]
Since $j_\psi(\nu)\geq0$, we obtain $j_\psi(\nu)=0$ on $S$ and hence $\psi_{|_S}=0$ by positive-definiteness of the Hermitian inner product $(\varphi,\phi)\mapsto\langle\nu\cdot\varphi,\phi\rangle$.
Since $\psi$ is uniquely determined by its values along a Cauchy hypersurface, we obtain $\psi=0$ on $M$.
\findemo
\end{proof}

Proposition \ref{obstruction} implies that if the initial data allow for a conformal extension and are not pure Maxwell theory, then the system has vanishing total charge.

\section{Proof of the main theorem}\label{proofmain}

%$y (t):= \frac{e^{t-c}}{1- e^{t-c}}$. The positive-valued branch of the function is its restriction to $(- \infty, c)$, and, as $y(c-a) > e^{-a}$, the distance to its pole can be estimated from below by the logarithm of the initial value. So, if $\Phi$ is semilinear and with punctured nonlinearity, then the lifetime $T_0 $ grows logarithmically with the inverse of the $H^k$ norm of $f$. \hfill \qed

%\section{Consequences for Dirac-Maxwell Theory}

\v In a first geometric step, we choose a $C^k$ extension $F$ of $(I^+(S'),g)$ to a globally hyperbolic manifold $(N,h)$ and consider the chosen Cauchy surface $S \subset I^+(S') $. 
Note that $ U:=  N \setminus J^-(S) $ is a future subset of $N$ and thus globally hyperbolic; let us choose a Cauchy temporal function $T$ on $U$, and consider a sequence of Cauchy hypersurfaces $S_n := T^{-1} (r_n)$ of $(U,h)$. The exact values of the $r_i$ will be specified later. Note that the $S_n$ are never Cauchy hypersurfaces of $F(I^+(S'))$. In the following we adopt the convention of denoting different spatio-temporal regularities explained after Theorem \ref{t:exsymmhypsys1storder}, related to the splitting induced by the temporal function $T$. The term $C^l H^k$ in this notation refers to an object which is $C^l$ regular in the time coordinate and $H^k$-regular in spatial direction. 

\v %Then we use the generalization of the $C^1$ extension criterion in Taylor for $C^k$ coefficients to show that there is a $C^1$ solution in $N$ and thus in $M$, as follows: 
The general strategy in the following is to find appropriate bounds on the initial values in different subsets of $F(S)$ (or, equivalently, corresponding bounds on $S$) implying that there is a global solution of a certain regularity. In our main theorem, we assume the initial Lorenz gauge condition on $F(S)$ (see Proposition \ref{p:constrLorenz}) and therefore can use the first prolongation (for the definition, see end of Appendix, after Corollary \ref{c:lifetime}) $\tilde{P}_{{\rm DW}}$ of the Dirac-wave operator $P_{{\rm DW}}$ in $N$ instead of $P_{{\rm DM}}$. We are first interested in regularity $C^1 H^4 $, as the degree of the operator $P_{{\rm DW}}$ is $2$ and as the critical regularity of the associated symmetric hyperbolic operator defined as a first prolongation is $k=4$ satisfying $\frac{k-1}{2} = 3/2$. 
Due to the lifetime estimate in Theorem \ref{t:lifetime}, which is a generalization of the well-known extension/breakdown criterion for smooth coefficients, there is a positive number $\d$ such that for initial values $u_1$ on $S_1$ with $\| u_1 \|_{H^s(S_1,h)} < \d$ there is a global solution on $D^+(S_1) \cap F(I^+(S))$ in $N$. Now, in a second step, we have to manage the ``initial jump'' from $S_\infty := F(S) $ to $ S_1 $, that is, we have to define sufficient conditions on $S_\infty$ such that initial values satisfying those conditions induce solutions $u$ reaching $S_1$ and satisfying $\|  u \|_{H^s(S_1,h)} < \d$ there, so we get a global solution on $D^+(F(S))$, where $D^+$ is the future domain of dependence. In the end, via conformally back-transforming the solution, we will obtain a solution on $J^+(S)$ with the given initial values on $S$.\\

Due to the unavoidable divergence of the conformal structure, we have to ``avoid spatial infinity'' in all computations, in the following sense: We transport sufficient $H^4$ bounds from $S_1$ down to $S_{\infty}  $ in regions of a certain distance from the boundary of $D^+(F(S)) \subset N$, while closer to the boundary we only transport them ``halfway down'' from one hypersurface $S_{n}$ to the next hypersurface $S_{n+1}$. More exactly, we choose a compact exhaustion of $S_{\infty}$, i.e. a sequence of open sets $C_n $ in $S_{\infty}$ such that $\overline{C_n}$ is compact, such that $\overline{C_n} \subset C_{n+1} $ and $\bigcup_{i=1}^{\infty} C_i =S_\infty$. Furthermore, we define $K_n := D^+ (C_n)$ as their future domains of dependence. We choose $ r_1 < \sup ( T(D^+(C_1) )) $. Inductively, by compactness of the possibly empty subset 

$$V_n:= J^+ (\overline{C_n}) \cap \partial K_{n+1}  ,$$ 

we find $\tau_n := {\rm min} \{ T(x) \vert x \in V_n \}  > - \infty $ and define $r_{n+1} :=   {\rm min} \{ r_n -1, \tau_n \} $ and $S_{n+1} := T^{-1} (r_{n+1}) $. With this choice, $\lim_{n \rightarrow \infty} r_n =  - \infty$ and

\begin{eqnarray}
\label{AbstandBitte!}
 J^- (S_{n+1} \setminus K_{n+1}) \cap C_n = \emptyset  .
\end{eqnarray}

Now we construct inductively a locally finite family of subsets $A_j $ of $F(S) $ and a sequence $b$ such that if $u_\infty$ is an initial value on $S_\infty$ with $\| u_\infty \vert_{A_j} \|_{C^4} < b_j $ then there is a global $C^1$ solution $u$ on $D^+(S_\infty) $ of $\tilde{P}_{DW}u = 0$ with $u \vert_{S_\infty} = u_\infty$. This sequence $b$ will be constructed via a corresponding sequence $a$ for the $H^4$ norms, which in turn is constructed as a limit of finite sequences $a^{(m)} \in \R^{m+1}$ that are stable in the sense that $a_n^{(m)} = a_n^{(m')}$ whenever $n \leq m-2, m'-2$, so that, for $n$ fixed, the sequence $m \mapsto a_{n}^{(m)}$ is eventually constant, thus we will,
 indeed, be able to define $a_i := \lim_{m \rightarrow \infty} a_i^{(m)}$ which will be a positive sequence.
 
 \bigskip
 
We define, for $n \geq1$, a finite set of subsets $\{ A_1^{(n)} , ... A_{n+1}^{(n)} \}  $ of $D^+(F(S))$ by (see figure below)

\bean
A_1^{(n)} := C_1 , \qquad A_{i+1}^{(n)} := J^- (S_i \setminus K_i) \cap  C_{i+1} \qquad \forall 1 \leq i \leq n-1 , \qquad A_{n+1}^{(n)} := J^-(S_n \setminus K_n) \cap S_{n+1} . 
\eean 

\begin{figure}[h]\label{fig:inductiveseq}
\begin{center}
\begin{tikzpicture}
%\caption{figure blabla};
\clip (-5,-1) rectangle (6,6);
\draw[very thick] (-5,0) -- (5,0);
\draw (5,0) node[anchor=west] {$S_\infty$};
\draw[very thick] (-5,0) -- (0,5);
\draw[very thick] (0,5) -- (5,0);
\draw (0.1,5) node[anchor=west] {$D^+(S_\infty)=F(I^+(S))$};
\draw[very thick] (-2,0) -- (0,2);
\draw (0,2) node[anchor=south] {$K_1$};
\draw[very thick] (0,2) -- (2,0);
\usetikzlibrary{snakes}\draw[snake=brace, raise snake=2pt, blue] (-2,0) -- (2,0);
\draw (0,0) node[anchor=south,color=blue] {\scriptsize$C_1= A_1^{(1)}$};
\draw[very thick] (-4,0) -- (0,4);
\draw (0,4) node[anchor=south] {$K_2$};
\draw[very thick] (0,4) -- (4,0);
\usetikzlibrary{snakes}\draw[snake=brace, mirror snake, raise snake=2pt] (-4,0) -- (4,0);
\draw (0,-0.1) node[anchor=north] {$C_2$};
\draw[very thick] (-5,1.2) -- (5,1.2);
\draw (5,1.2) node[anchor=west] {$S_1$};
\draw[thick,style=dashed] (-2,0) -- (-3,1);
\draw[thick,style=dashed] (2,0) -- (3,1);
\draw[very thick] (-5,0.7) -- (5,0.7);
\draw (5,0.7) node[anchor=west] {$S_2$};
\draw[thick,style=dashed] (-0.8,1.2) -- (-0.3,0.7);
\draw[thick,style=dashed] (0.8,1.2) -- (0.3,0.7);
\draw[very thick,color=blue] (-0.3,0.7) -- (-4.3,0.7);
\draw[very thick,color=blue] (0.3,0.7) -- (4.3,0.7);
\usetikzlibrary{snakes}\draw[snake=brace, raise snake=1pt, blue] (0.3,0.7) -- (4.3,0.7);
\draw (2.3,0.7) node[anchor=south,color=blue] {\scriptsize$A_2^{(1)}$}; 
\usetikzlibrary{snakes}\draw[snake=brace, raise snake=1pt, blue] (-4.3,0.7) -- (-0.3,0.7);
\draw (-2.3,0.7) node[anchor=south,color=blue] {\scriptsize$A_2^{(1)}$}; 
\end{tikzpicture}
\caption{Construction of the sequence $A_i^{(n)}$}
\end{center}
\end{figure}
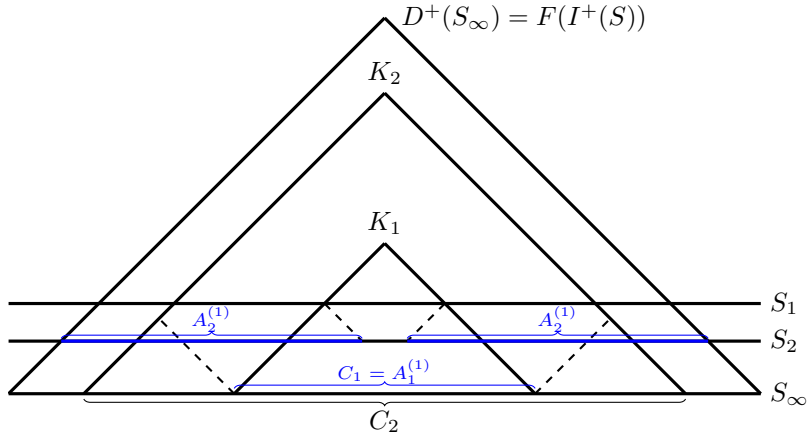
   
Note that the first $n$ subsets are in $S_\infty$ whereas the last one is in $S_{n+1}$. Note furthermore that the sequence stabilizes in the sense that $A_i^{(n)} = A_i^{(m)}$ if $m,n > i+1$, and the limit sequence is $A_1 := C_1,  A_{i+1} := J^- (S_i \setminus K_i) \cap C_{i+1} \qquad \forall i>1$.

Let us call a finite positive sequence $a_1^{(n)} , ... , a_{n+1}^{(n)}$ a {\bf control sequence at step $n$} iff every $C^1 H^4$ solution $u$ of $\tilde{P}_{DW}u=0$ in
$J^+(S_\infty) \cap J^-(S_{n+1})$ with $\| u \vert_{A_i^{(n)}} \|_{H^4} < a_i^{(n)}$ for all $i \in \N \cap [0,n+1]$ extends to a global $C^1 H^4$ solution on $D^+(S_\infty) = F(I^+(S))$.
 
\begin{Lemma} 
For every $ n \geq 2 $, there is a control sequence $a^{(n)}_i$ at step $n$, and the sequences stabilize in the sense that $a^{(n)}_i = a^{(m)}_i$ if $m,n > i+1  $.
\end{Lemma}

 \V {\bf Proof of the lemma.} Obviously, for $n=1$, we only have to ensure that $ \vert \vert u \vert \vert_{H^4 (S_1)} \leq \delta$. The lifetime estimate of Theorem \ref{t:lifetime} in the region $K_1$ implies that there is a positive constant $a_1^{(2)} $ such that $ \vert \vert u \vert \vert_{H^4(C_1)} < a_1^{(2)}$ ensures that $u$ extends up to $S_1 \cap K_1$ and $\vert \vert u \vert \vert_{H^4 (S_1 \cap K_1)} < \delta/2 $. Moreover, the lifetime estimate in $I^+(S_2) \cap I^-(S_1 \setminus K_1)$ implies that there is a second constant $a_2^{(2)}$ such that $ \vert \vert u \vert \vert_{ H^4(S_2 \cap J^-(S_1 \setminus K_1)) }  < a_2^{(2)}$ implies $  \vert \vert u \vert \vert_{H^4(S_1 \setminus K_1)} < \delta/2  $. Then it is straightforward to show that if both conditions are satisfied, the solution $u$ fulfills $\vert \vert u  \vert \vert_{H^4 (S_1)} < \delta$, and therefore the solution extends to all of $F(J^+(S))$.
 
 \bigskip
 
Each induction step is again done by applying the lifetime estimate in two regions. Now assume that there is a control sequence at step $n$. We have to look for an appropriate sequence of $H^4$ bounds $(a^{(n+1)}_1, .. , a^{(n+2)}_{n+2})$ on $A_{1}^{(n+1)}, ... , A_{n+2}^{(n+1)}$. First we define

\bean
a_i^{(n+1)} := a_i^{(n)} \qquad \forall 1 \leq i \leq n
\eean

\bigskip

To ensure the $H^4$-bound on $A_{n+1}^{(n)}$, we divide $A_{n+1}^{(n)}$ into its inner part $I_{n+1}^{(n)} := A_{n+1}^{(n)} \cap K_{n+1} $ and its outer part $O_{n+1}^{(n)} := A_{n+1}^{(n)} \setminus K_{n+1} = S_{n+1} \setminus K_{n+1}$. We want to ensure the $H^4$-bound $a_{n+1}^{(n)}$ on both parts. To guarantee the $H^4$ bound $a^{(n)}_{n+1}$ on the inner part there is a sufficient $H^4$ bound $a_{n+1}^{(n+1)}$ on 

\bean
J^-(I_{n+1}^{(n)}) \cap S_\infty = J^-(S_n \setminus K_n) \cap C_{n+1} = A_{n+1}^{(n+1)} ,
\eean

whereas for the $H^4$ bound $a^{(n)}_{n+1}$ on the outer part, an $H^4$ bound $ a_{n+2}^{(n+1)} $ on 

\bean
J^-(O_{n+1}^{(n)}) \cap S_{n+2} = J^-(S_{n+1} \setminus K_{n+1}) \cap S^+_{n+2} = A_{n+2}^{(n+1)}
\eean

is sufficient. Thus $a_i^{(n+1)}$ is a control sequence at step $n+1$, and indeed the sequences stabilize in the sense above by definition. \hfill \qed

\bigskip

\V As the sequences stabilize, we can define the (infinite, positive) limit sequence $a_i$. Now there are $b_i>0$ such that $\| u_0 \| _{H^4 (A_i)} < a_i $ is satisfied if $\| u_0 \| _{C^4 (A_i)} < b_i $. Now, the condition \ref{AbstandBitte!} ensures that for the annular regions $D_i:= C_{i+1} \setminus C_i$, with $D_0:= C_0$ and for every $i \in \N$ we have $D_i \cap A_j \neq \emptyset$ only if $j=i $ or $j=i+1$. So on every $D_i$ we have to satisfy only two $C^4$ bounds $b_i$ for all control sequences to be satisfied; let $\underline{b}_i$ be the minimum of those two bounds. Now, given initial values $u_0$ with 
\bea
\label{weighted}
 \vert \vert u_0  \vert \vert_{C^4(D_i)} < \underline{b}_i , 
\eea
and given any point $q \in F(M)$, we want to show that $q$ is contained in a domain of definition for a $C^1$ solution $u$ of $\tilde{P}_{DW}u=0$ with $u \vert_{S_\infty} = u_0$. To that purpose, we choose an $i$ such that $q \in K_i$ and choose $f_i \in C^{\infty} (S_\infty, [0,1])$ with $f_i (C_i ) = \{ 1 \} $ and ${\rm supp} (f_i ) \subset S_\infty \setminus C_{i+1}$. Then we solve the initial value problem for $u^{(i)} = f_i \cdot u_0$. Applying the $i$th step 
in the induction above, we get a solution $u^{[i]}$ on a domain of definition including $q$. Locality implies that any local solution with initial value $u_0$ coincides with $u^{[i]}$ on $K_i$. This is, the domain of definition of a maximal solution includes $q$. Note that Eq. \ref{weighted} corresponds to a bound in a weighted $C^4$-space on $S$.

%Let $S$ be the image in $E$ of a Cauchy surface in $M$. For given $\e >0$, there is a positive number $\d$ such that for initial values $a$ on $B_{\e}^+ (S)$ with $\vert \vert a \vert \vert_{H^s} < \d$ there is a global solution on $M$. Now, in a second step, we have to manage the 'initial jump' from $S$ to $B_{\e}^+ (S)$. This is done via calculations in $M$: For every  $ p \in B^+_{\e} (S) \cap c(M)$ is the image under $C$ we have $I_E^- (p) = c (I_M^- (q))$ for $c(q) =p$, but the size of the past cones grows to infinity towards $i_0$. Therefore we have to show that the existence time grows to infinity, with the inverse power of the conformal factor. This is done by a theorem on the growth of existence time for equations with quadratic nonlinearity. The necessity of the operator being cosp is easily seen by the ODE $x' = 1 + x^2$ solved by $tan$ which blows up in time $\pi$ at the latest irrespectively of the initial value.

\bigskip

\v As usual, we show higher regularity by bootstrapping, i.e. considering the differentiated equation (which is a linear equation in the highest derivatives again). Consider the highest derivatives in a Sobolev Hilbert space as independent variables and show that they are in the same Sobolev Hilbert space as the coefficients, thereby gaining one order of (weak) differentiability. Finally we use Sobolev embeddings in the usual way.  \hfill\qed

\bigskip

\section{Appendix: Modification of the breakdown criterion, existence time and regularity}\label{s:appendix}

Following \cite[Ch. 16]{Taylor3} but modifying the proof so as to allow for coefficients of finite regularity, we present the proof of local existence and uniqueness for solutions to symmetric hyperbolic systems.
Although we could not find the existence theory for symmetric hyperbolic systems with coefficients of finite regularity in the literature, we do not claim originality of the following results but present them in full detail for the sake of self-containedness.

\bdefi[\protect{\cite[Sec. 16.2]{Taylor3}}]\label{d:symmhypsys1storder}
For $\mathbb{K}=\R$ or $\C$ and $N\in\mathbb{N}$, a first-order \emph{\bf symmetric hyperbolic system on $\R^n$  with values in $\mathbb{K}^N$} is a system of equations of the form
\begin{equation}\label{eq:symmhypsys1storder}
\left\{\begin{array}{ll}A_0(t,x,u)\frac{\partial u}{\partial t}&=L(t,x,u,\partial)u+g(t,x,u)\textrm{ on }\R\times\R^n\\ u(0)&=f,\end{array}\right. 
\end{equation}
where 
\begin{itemize}\item $L(t,x,u,\partial)v:=\sum_{j=1}^nA_j(t,x,u)\partial_jv$ for all $v:\R\times\R^n\to\mathbb{K}^N$, with $A_j:\R\times\R^n\times\mathbb{K}^N\to\mathrm{Mat}_{N\times N}(\mathbb{K})$ such that $A_j^*=A_j$ (pointwise),
\item $A_0:\R\times\R^n\times\mathbb{K}^N\to\mathrm{Mat}_{N\times N}(\mathbb{K})$ such that $A_0^*=A_0$ (pointwise) and $A_0(t,x,u)\geq c\cdot\mathrm{I}$ for some $c>0$,
\item $g:\R\times\R^n\times\mathbb{K}^N\to\mathbb{K}^N$ and
\item $f:\R^n\to\mathbb{K}^N$.
\end{itemize}
\edefi

The same definition can be made when replacing $\R^n$ by an $n$-dimensional torus $\mathbb{T}^n$.
The condition on $A_0$ means that $A_0$ is a pointwise Hermitian/symmetric matrix that is \emph{uniformly} positive definite on $\R\times\R^n\times\mathbb{K}^N$.\\

We want to prove the local existence and the uniqueness of solutions to first-order symmetric hyperbolic systems on $\mathbb{T}^n$.
Later on, we shall consider the case of higher order symmetric hyperbolic system also on other manifolds.\\
%\footnote{To do, in particular on $\R^n$.}.\\

We start by assuming low regularity on the data (we shall see below how the regularity of the solution depends on that of the data).
The main theorem we want to prove is the following:

\btheo\label{t:exsymmhypsys1storder}
Consider a $\mathbb{K}^N$-valued first-order symmetric hyperbolic system on $\mathbb{T}^n$ as in {\rm Definition \ref{d:symmhypsys1storder}} and assume $A_j,g$ to be $C^1$ in $(t,x,u)$ and $C^k$ in $(x,u)$ for some $k>\frac{n}{2}+1$.
%\footnote{Need probably $C^1$ in $(t,x,u)$ and even $C^k$ in $(x,u)$.} 
Then for any $f\in H^k$,
\begin{enumerate}
\item[$1.$] there is an $\eta\in\R_+^\times$ for which a unique solution $u\in C^1(]-\eta,\eta[\times\mathbb{T}^n)\cap C^0(]-\eta,\eta[,H^k(\mathbb{T}^n))$ to {\rm (\ref{eq:symmhypsys1storder})} exists;
%\footnote{Kann doch nur in linsenf\"ormigen bzw. Abh\"angigkeitsgebieten eindeutig sein!}  
\item[$2.$] {\em (extension criterion):} that solution $u$ exists as long as $\|u(t)\|_{C^1(\mathbb{T}^n)}$ remains bounded.
\end{enumerate}
\etheo

By $C^0H^k$, we mean continuous in the first variable $t\in I$ with values in the $H^k$-Sobolev space on $\mathbb{T}^n$ or $\mathbb{T}^n\times\mathbb{K}^N$.
We shall mostly omit the interval $I$ or the torus $\mathbb{T}^n$ in the notation.
As usual, $H^k:=W^{k,2}$.
In the sequel, we shall often denote those spaces of functions with regularity $R$ in $t$ and with values in a Banach space $S$ (mostly of functions in the other variables) with $RS$ (e.g. $C^0H^k$, $L^\infty H^k$ etc.).\\

\bigskip

During the seven-step proof of Theorem \ref{t:exsymmhypsys1storder}, in several estimates, as multiplicative factors functions $C_i: \R^m \rightarrow \R$ will appear that take certain norms of the (approximate) solutions or of other maps as arguments. For simplicity, we will adopt the convention that these functions ('constants only depending on the norm') are taken to be {\em monotonously increasing}, and we try to number them consecutively by indices in every of the seven steps of the proof, which are the following:

\begin{enumerate}
\item Using mollifiers, perturb (\ref{eq:symmhypsys1storder}) by a small parameter $\varepsilon>0$ in order to obtain a new system that can be interpreted as an ODE in the Banach space $H^k=H^k(\mathbb{T}^n)$.
\item For each value of the parameter $\varepsilon>0$, solve the corresponding ODE locally about $0\in\R$ and obtain a so-called approximate solution.
\item By a uniform (in the parameter $\varepsilon$) control of the pointwise $H^k$-norm of those approximate solutions, show that they all exist on a common interval $]-\eta,\eta[$ with $\eta>0$.
\item Up to shrinking $\eta$ a bit, extract of the families of approximate solutions a weak accumulation point and show that it is a $C^1$-solution to (\ref{eq:symmhypsys1storder}) on $]-\eta,\eta[\times\mathbb{T}^n$.
\item Show uniqueness of the local solution by controlling the rate of convergence of the approximate solutions against the solution when $\varepsilon\to0$.
\item Improve the regularity of the solution to $C^0H^k$.
This proves $1.$
\item Show that in fact $\|u(t)\|_{H^k(\mathbb{T}^n)}$ remains bounded as long as $\|u(t)\|_{C^1(\mathbb{T}^n)}$ does. 
Assuming the solution $u$ stops existing at $T>0$, use a precise control of the length of the existence interval in the theorem of Picard-Lindel\"of to prove that all approximate solutions - for an initial value fixed ``shortly before'' $T$ - can be extended beyond $T$; this also implies (using uniqueness) that the solution can be extended beyond $T$, contradiction.
\end{enumerate}

Let $J_\varepsilon$ be the convolution with $\theta_\varepsilon=\varepsilon^{-n}\theta(\frac{\cdot}{\varepsilon})$, where $\theta\in C^\infty(\mathbb{R}^n,[0,\infty[)$, $\mathrm{supp}(\theta)\subset \overline{B}_1(0)$, $\int_{\R^n}\theta dx=1$ and $\theta\circ(-\mathrm{Id})=\theta$; the last condition is needed for the self-adjointness of $J_\varepsilon$ in $L^2$ and higher Sobolev spaces.
The operator $J_\varepsilon$ is a smoothing operator approximating the identity in the following sense: $J_\varepsilon\buil{\longrightarrow}{\varepsilon\searrow 0}\mathrm{Id}$ pointwise in $W^{k,q}(\mathbb{R}^n)$ for every $(k,q)\in\mathbb{N}\times[1,\infty[$ and also pointwise in $C^0(I,C^k(\mathbb{T}^n))$ for any open interval $I$.
We shall often make use of $[J_\varepsilon,\partial^\alpha]=0$ for every multi-index $\alpha$ and of the following facts: $J_\varepsilon\colon W^{k,q}(\mathbb{R}^n)\to W^{k,q}(\mathbb{R}^n)\cap C^\infty(\mathbb{R}^n)$ has norm $\|J_\varepsilon\|\leq1$, $J_\varepsilon\colon C_b^k(\mathbb{R}^n)\to C_b^k(\mathbb{R}^n)$ has norm  $\|J_\varepsilon\|\leq1$, the operator $J_\varepsilon\colon C_b^0(\mathbb{R}^n)\to C_b^k(\mathbb{R}^n)$ has norm $\|J_\varepsilon\|\leq C(\varepsilon)$,  the operator $J_\varepsilon\colon C^0(I,C^k(\mathbb{T}^n))\to C^0(I,C^k(\mathbb{T}^n))$ has norm $\|J_\varepsilon\|\leq 1$.
It is also interesting to notice that $J_\varepsilon$ is an operator $\mathrm{Lip}(I,H^{k-1})\to\mathrm{Lip}(I,H^k)$ with $\|J_\varepsilon u\|_{C^{0,1}(I,H^k)}\leq C(\varepsilon)\|u\|_{C^{0,1}(I,H^{k-1})}$, where $I\subset\mathbb{R}$ is a bounded open interval and $\mathrm{Lip}(I,H^l)=C^{0,1}(I,H^l)$.
Namely for any $f\in H^{k-1}$ and $\alpha\in\mathbb{N}^n$ with $|\alpha|\leq k$, one has
\begin{eqnarray*} \|\partial^\alpha J_\varepsilon f\|_{L^2}^2&=&\|f*\partial^\alpha\theta_\varepsilon\|_{L^2}^2\\
&\leq&\|f\|_{L^2}^2\cdot\underbrace{\|\partial^\alpha\theta_\varepsilon\|_{L^1}^2}_{C(\varepsilon)^2}\\
&\leq&C(\varepsilon)^2\|f\|_{H^{k-1}}^2,
\end{eqnarray*} 
so that $\|J_\varepsilon f\|_{H^k}\leq C(\varepsilon)\|f\|_{H^{k-1}}$, which shows the claim.

\bprop\label{p:IminusJeps}
$\|\mathrm{Id}-J_\varepsilon\|_{\mathcal{L}(H^{1},L^2)}\leq C\cdot\varepsilon$.\footnote{The statement holds as well for $\R^n$ instead of $\mathbb{T}^n$ with the same proof \emph{mutatis mutandis}.}
\eprop
%\footnote{Taylor claims \cite[eq. (1.38)]{Taylor3} that in fact $\|\mathrm{Id}-J_\varepsilon\|_{\mathcal{L}(H^{k},L^2)}\leq C\cdot\varepsilon^k$, but I momentarily cannot prove it (N).}\\
\begin{proof} For any $f\in H^1(\mathbb{T}^n)$, we have 
\begin{eqnarray*} 
\|J_\varepsilon f-f\|_{L^2}^2&\leq&\int_{\mathbb{T}^n}\left(\int_{\mathbb{T}^n}|f(x-y)-f(x)|\cdot\theta_\varepsilon(y)dy\right)^2dx\\
&\leq&\int_{\mathbb{T}^n}\left(\int_{B_\varepsilon(0)}|f(x-y)-f(x)|^2dy\right)\|\theta_\varepsilon\|_{L^2(B_\varepsilon)}^2dx\\
&\leq&\varepsilon^{-n}\|\theta\|_{L^2}^2\cdot\int_{\mathbb{T}^n}\left(\int_{B_\varepsilon(0)}|f(x-y)-f(x)|^2dy\right)dx\\
&\leq&\varepsilon^{-n}\|\theta\|_{L^2}^2\cdot\int_{\mathbb{T}^n}\int_{B_\varepsilon(0)}\left(\int_0^1|d_{x-ty}f(y)|^2dt\right)dydx\\
&\leq&\varepsilon^2\cdot\varepsilon^{-n}\|\theta\|_{L^2}^2\cdot\int_{\mathbb{T}^n}\int_{B_\varepsilon(0)}\left(\int_0^1|d_{x-ty}f|^2dt\right)dydx\\
&\leq&\varepsilon^{2-n}\|\theta\|_{L^2}^2\cdot\int_{B_\varepsilon(0)}\int_0^1\|df\|_{L^2(\mathbb{T}^n)}^2dtdy\qquad\textrm{(Fubini)}\\
&\leq&\varepsilon^{2-n}\|\theta\|_{L^2}^2\cdot\mathrm{Vol}(B_\varepsilon(0))\cdot\|df\|_{L^2(\mathbb{T}^n)}^2\\
&\leq&C\cdot\|\theta\|_{L^2}^2\cdot\varepsilon^2\cdot\|df\|_{L^2(\mathbb{T}^n)}^2\\
&\leq&C\cdot\varepsilon^2\cdot\|f\|_{H^1(\mathbb{T}^n)}^2,
\end{eqnarray*} 
which concludes the proof of the proposition.
\hfill\qed\\
\end{proof}

In the proof of Theorem \ref{t:exsymmhypsys1storder}, we use the following inequalities, see e.g. \cite[Prop. 13.3.7]{Taylor3}, \cite[Thm. 2.2.2, 2.2.3 \& Lemma 2.2.6]{Rendall06} and \cite[Thm. 2.3.6 \& 2.3.7]{Finster02}.

\blemm[Moser]\label{l:Moserestim}
Let $k,n\in\mathbb{N}\setminus\{0\}$.
\begin{itemize}
\item[i)] {\bf (First Moser estimate)} There exists a constant $C=C(k,n)\in\R_+^\times$ such that, for all $f,g\in L^\infty(\R^n)\cap H^k(\R^n)$,
\begin{equation}\label{eq:Moser1} 
\|f\cdot g\|_{H^k}\leq C\cdot\left(\|f\|_{L^\infty}\|g\|_{H^k}+\|f\|_{H^k}\|g\|_{L^\infty}\right). 
\end{equation}
\item[ii)] {\bf (Second Moser estimate)} There exists a constant $C=C(k,n)\in\R_+^\times$ such that, for all $f\in W^{1,\infty}(\R^n)\cap H^k(\R^n)$, $g\in L^\infty(\R^n)\cap H^{k-1}(\R^n)$ and $\alpha\in\mathbb{N}^n$ with $|\alpha|\leq k$,
\begin{equation}\label{eq:Moser2}
\|\partial^\alpha(fg)-f\partial^\alpha g\|_{L^2}\leq C\cdot\left(\|\nabla f\|_{H^{k-1}}\|g\|_{L^\infty}+\|\nabla f\|_{L^\infty}\|g\|_{H^{k-1}}\right).
\end{equation}
\item[iii)] {\bf (Third Moser estimate)} Let $F\in C^\infty(\mathbb{K}^N,\mathbb{K}^L)$ with $F(0)=0$.
Then there is a constant $C\in\R_+^\times$, which only depends on $k,n,F$ and on $\|f\|_{L^\infty}$, such that, for any $f\in L^\infty(\R^n)\cap H^k(\R^n)$ and $\alpha\in\mathbb{N}^n$ with $|\alpha|\leq k$,
\begin{equation}\label{eq:Moser3}
\|\partial^\alpha F(f)\|_{L^2}\leq C(\|f\|_{L^\infty})\cdot\|\nabla^{|\alpha|} f\|_{L^2}.
\end{equation}
\end{itemize}
\elemm 

In \cite[Prop. 13.3.9]{Taylor3}, there is the following alternative (and weaker) version of (\ref{eq:Moser3}): for every $F\in C^\infty(\mathbb{K}^N,\mathbb{K}^L)$ with $F(0)=0$, there exists a constant $C>0$ depending only on $k,n,F$ and on $\|f\|_{L^\infty}$ such that, for all $f\in L^\infty(\R^n)\cap H^k(\R^n)$,
\begin{equation}\label{eq:Moser3bis} 
\|F(f)\|_{H^k}\leq C(\|f\|_{L^\infty})\cdot(1+\|f\|_{H^k}). 
\end{equation}

Note that all estimates from Lemma \ref{l:Moserestim} remain true when replacing $\R^n$ by the $n$-dimensional torus $\mathbb{T}^n$.
Moreover, since $\mathbb{T}^n$ has finite volume, the assumption $F(0)=0$ can be dropped for the weaker third Moser estimate (\ref{eq:Moser3bis}), however not for (\ref{eq:Moser3}) and $\alpha=0$.

\blemm\label{l:mollifiercommutator}
Let $A\in C^1(\R^n)$, $p\in[1,\infty[$ and $\varepsilon>0$.
Then there exists a constant $C=C(n,p)>0$ such that, for any $v\in L^p(\R^n)$,
\begin{itemize}
\item[i)] $\|[A,J_\varepsilon]v\|_{L^p}\leq\left\{\begin{array}{c}\vspace{.2cm}C\cdot\|A\|_{C^0}\cdot\|v\|_{L^p}\\C\cdot\varepsilon\cdot\|A\|_{C^1}\cdot\|v\|_{L^p}\end{array}\right.$ .
\item[ii)] $\|[A,J_\varepsilon]v\|_{W^{1,p}}\leq C\cdot\|A\|_{C^1}\cdot\|v\|_{L^p}$.
\item[iii)] $\|[A,J_\varepsilon]\frac{\partial v}{\partial x_j}\|_{L^p}\leq C\cdot\|A\|_{C^1}\cdot\|v\|_{L^p}$.
\end{itemize}
\elemm 

\begin{proof} See e.g. \cite[Ex. 13.1.1 - 13.1.3]{Taylor3}.
\qed\\
\end{proof}

%{\red\bf domain of definition? Nicolas}.\\

{\bf Step 1}: We mollify the symmetric hyperbolic system in order to obtain an ODE in $H^k$.\\
% We have to write the corresponding part in diracmaxwell4 down because of the new factor $A_0$.
% Check for whether $H^k$ regularity of the data suffices or not; after discussions, we think we probably need $A_j$ to be $C^k$ in $(x,u)$ and also $g$  $C^k$ in $(x,u)$.\\

{\bf Claim 1}: {\em For any sufficiently small $\varepsilon>0$, the equation $A_0(t,x,J_\varepsilon u_\varepsilon)\frac{\partial u_\varepsilon}{\partial t}=J_\varepsilon L(t,x,J_\varepsilon u_\varepsilon)J_\varepsilon u_\varepsilon+J_\varepsilon g(t,x,J_\varepsilon u_\varepsilon)$ is an ODE in $H^k$ that is strongly locally Lipschitz in $u_\varepsilon$, that is, there exists a Lipschitz constant (in $x$) on all products $[0,T]\times \overline{B}_R(0)$, where $\overline{B}_R(0)$ is the closed $R$-ball about $0\in H^k$.}\\

\begin{proof} Consider the map $F\colon\R\times H^k\to H^k$, 
\[F(t,v)(x):=A_0^{-1}(t,x,(J_\varepsilon v)(x))\cdot\Big(J_\varepsilon\big[y\mapsto L(t,y,(J_\varepsilon v)(y))(J_\varepsilon v)(y)\big](x)+J_\varepsilon\big[y\mapsto g(t,y,(J_\varepsilon v)(y))\big](x)\Big)\]
for all $(t,v)\in\R\times H^k$ and every $x\in\mathbb{T}^n$.
As in \cite{Taylor3}, we shortly write
\[F(t,v)=A_0^{-1}(t,x,J_\varepsilon v)\cdot\left(J_\varepsilon L(t,x,J_\varepsilon v)J_\varepsilon v+J_\varepsilon g(t,x,J_\varepsilon v)\right)\]
for every $v\in H^k$.
We show that $F$ is $C^1$ (in the Fr\'echet sense) with bounded differential on each subset of the form $[0,T]\times \overline{B}_R(0)$ in $I\times H^k$.
We only treat the case of one term in the definition of $F$, the others being handled in a similar manner.
Namely consider the map $(t,v)\mapsto J_\varepsilon g(t,x,J_\varepsilon v)$ from $\R\times H^k\to H^k$.
Then for any $v\in H^k$ and $h\in H^k$, we have 
\begin{eqnarray*} 
J_\varepsilon g(t,x,J_\varepsilon (v+h))-J_\varepsilon g(t,x,J_\varepsilon v)&=&J_\varepsilon\left(g(t,x,J_\varepsilon (v+h))-g(t,x,J_\varepsilon v)\right)\\
&=&J_\varepsilon\left(x\mapsto g_u'(t,x,(J_\varepsilon v)(x))\cdot(J_\varepsilon h)(x)+|(J_\varepsilon h)(x)|\cdot\epsilon((J_\varepsilon h)(x))\right),
\end{eqnarray*}
where $\epsilon(w)\buil{\longrightarrow}{w\to 0}0$.
The map $h\mapsto J_\varepsilon\left(x\mapsto g_u'(t,x,(J_\varepsilon v)(x))\cdot(J_\varepsilon h)(x)\right)$ is linear and  bounded $H^k\to H^k$:
\begin{eqnarray*} 
\left\|J_\varepsilon\left(x\mapsto g_u'(t,x,(J_\varepsilon v)(x))\cdot(J_\varepsilon h)(x)\right)\right\|_{H^k}&\leq&C_1(\varepsilon)\left\|x\mapsto g_u'(t,x,(J_\varepsilon v)(x))\cdot(J_\varepsilon h)(x)\right\|_{L^2}\\
&\leq&C_1(\varepsilon)\left\|x\mapsto g_u'(t,x,(J_\varepsilon v)(x))\right\|_{L^\infty}\cdot\left\|J_\varepsilon h\right\|_{L^2}\\
% &\bui{\leq}{\rm(\ref{eq:Moser1})}&C(\varepsilon)\Big(\left\|x\mapsto g_u'(t,x,(J_\varepsilon v)(x))\right\|_{H^{k-1}}\cdot\left\|J_\varepsilon h\right\|_{L^\infty}\\
% &&\phantom{C(\varepsilon)\Big(}+\left\|x\mapsto g_u'(t,x,(J_\varepsilon v)(x))\right\|_{L^\infty}\cdot\left\|J_\varepsilon h\right\|_{H^{k-1}}\Big)\\
&\leq&C_1(\varepsilon,t)\cdot\|h\|_{H^k},
\end{eqnarray*} 
where we have used the compactness of $\mathbb{T}^n$ and the fact that $g_u'$ is continuous.
Furthermore, the map $h\mapsto J_\varepsilon\left(x\mapsto |(J_\varepsilon h)(x)|\cdot\epsilon((J_\varepsilon h)(x))\right)$ is of the form $\mathrm{o}(\|h\|_{H^k})$ since 
\begin{eqnarray*} 
\frac{\left\|J_\varepsilon\left(x\mapsto |(J_\varepsilon h)(x)|\cdot\epsilon((J_\varepsilon h)(x))\right)\right\|_{H^k}}{\|h\|_{H^k}}&\leq&\frac{C_1(\varepsilon)}{\|h\|_{H^k}}\left(\left\||J_\varepsilon h|\cdot\epsilon((J_\varepsilon h))\right\|_{L^2}\right)\\
% &\leq&\frac{C(\varepsilon)}{\|h\|_{H^k}}\Big(\|J_\varepsilon h\|_{H^{k-1}}\cdot\|\epsilon((J_\varepsilon h))\|_{L^\infty}\\
% &&\phantom{\frac{C(\varepsilon)}{\|h\|_{H^k}}\Big(}+\|J_\varepsilon h\|_{L^\infty}\cdot\|\epsilon((J_\varepsilon h))\|_{H^{k-1}}\Big)\\
&\leq&\|\epsilon((J_\varepsilon h))\|_{L^\infty},
\end{eqnarray*} 
where $\|\epsilon((J_\varepsilon h))\|_{L^\infty}\leq C_2(\varepsilon,t)\|h\|_{L^1}\buil{\longrightarrow}{\|h\|_{H^k}\to0}0$ because of $\mathbb{T}^n$ being compact.
Finally, the map 
\[H^k\to\mathcal{B}(H^k,H^k),\qquad v\mapsto J_\varepsilon\left(x\mapsto g_u'(t,x,(J_\varepsilon v)(x))\cdot(J_\varepsilon \bullet)(x)\right)\]
is continuous and bounded on each ball in $H^k$: this follows from the same kind of estimates as above as well as the continuity of $g_u'$ on $\R\times\mathbb{T}^n\times\mathbb{K}^N$.
%{\red\bf(Is it enough or do you want more?) }
This shows the claim.
%\begin{eqnarray*} &=&J_\varepsilon\left(x\mapsto\int_0^1g_u'(t,x,(J_\varepsilon v)(x)+s(J_\varepsilon h)(x))(J_\varepsilon h)(x)\,ds\right).
% \end{eqnarray*} 
% The map $D(t,v,h)\colon$ depends obviously linearly of $h\in H^k$ and we have 
% \begin{eqnarray*} 
% \left\|D(t,v,h)\right\|_{H^k}&\leq&\left\|x\mapsto\int_0^1g_u'(t,x,(J_\varepsilon v)(x)+s(J_\varepsilon h)(x))(J_\varepsilon h)(x)\,ds\right\|_{H^k}\\
% &\leq&\int_0^1\left\|x\mapsto g_u'(t,x,(J_\varepsilon v)(x)+s(J_\varepsilon h)(x))(J_\varepsilon h)(x)\right\|_{H^k}\,ds\\
% &\leq&\left\|(s,x)\mapsto g_u'(t,x,(J_\varepsilon v)(x)+s(J_\varepsilon h)(x))\right\|_{L^\infty}\cdot\left\|J_\varepsilon h\right\|_{H^k}\\
% &\leq&C\left(\|J_\varepsilon v\|_{L^\infty},\|J_\varepsilon h\|_{L^\infty}\right)\cdot\|h\|_{H^k}\\
% &\leq&
% \end{eqnarray*} 
%{\red\bf to be written down.}
\qed\\ 
\end{proof}

{\bf Step 2}: This is mainly classical ODE theory, applicable as soon as the nonlinearity is continuous (in $(t,x)$) and locally Lipschitz (in the usual sense) in $x$.\\

{\bf Claim 2}: {\em For any $f\in H^k$ and any sufficiently small $\varepsilon>0$ , the system 
\begin{equation}\label{eq:sysepsilon} 
\left\{\begin{array}{ll}A_0(t,x,J_\varepsilon u_\varepsilon)\frac{\partial u_\varepsilon}{\partial t}&=J_\varepsilon L(t,x,J_\varepsilon u_\varepsilon)J_\varepsilon u_\varepsilon+J_\varepsilon g(t,x,J_\varepsilon u_\varepsilon)\\ u_\varepsilon(0)&=f\end{array}\right. 
\end{equation}
has a unique solution $u_\varepsilon\in C^1(]-\eta_\varepsilon,\eta_\varepsilon[,H^k)$ for some $\eta_\varepsilon>0$.}\\

\begin{proof}
straightforward consequence of the theorem of Picard-Lindel\"of.
\qed\\ 
\end{proof}

{\bf Step 3}: ``Standard estimates'' based on Moser(-Trudinger) estimates and on Bihari's inequality \cite{Bihari56}.\\
% Here one should also pay attention to the fact that, in order for the ``standard'' extension criterion (boundedness of the pointwise norm of the solution) for ODEs to be applied, the nonlinearity should be more than locally Lipschitz: it should be locally Lipschitz in the sense of Step 1.
% Note that closed balls in $H^k$ are not compact, therefore this Lipschitz condition is not guaranteed by the usual local Lipschitz one.
% The stronger condition is fulfilled here by Step 1.\\

{\bf Claim 3}: {\em Under the assumptions of {\rm Claim 2} and with $k>\frac{n}{2}+1$, there exists an $\eta>0$ and a $K\in[0,\infty[$ such that $\|u_\varepsilon(t)\|_{H^k}\leq K$ for all $t\in]-\eta,\eta[$. 
In particular, the number $\eta_\varepsilon$ from {\rm Claim 2} may be chosen independently on $\varepsilon$.}\\

\begin{proof}
We introduce the new $L^2$-Hermitian inner product $(\cdot\,,\cdot)_{L^2,\varepsilon}:=(A_{0\varepsilon}\cdot\,,\cdot)_{L^2}$ on $\mathbb{T}^n$, where $A_{0\varepsilon}:=A_0(t,x,J_\varepsilon u_\varepsilon)$.
Note that $(\cdot\,,\cdot)_{L^2,\varepsilon}$ depends on $\varepsilon>0$ and also implicitely on $t$; but by assumption on $A_0$ and because we only consider compact sets of the form $[0,T]\times\mathbb{T}^n$, the norms $\|\cdot\|_{L^2,\varepsilon}$ and $\|\cdot\|_{L^2}$ are equivalent; more precisely, for any $T\in[0,\infty[$, there exists $C=C(T,\|u_\varepsilon\|_{C^0([0,T],L^\infty)})\in]0,\infty[$ such that $c\|\cdot\|_{L^2}^2\leq\|\cdot\|_{L^2,\varepsilon}^2\leq C\|\cdot\|_{L^2}^2$, where $c>0$ is the constant from Definition \ref{d:symmhypsys1storder}.
We pick an arbitrary $\alpha\in\mathbb{N}^n$ with $|\alpha|\leq k$ and estimate $\|\partial^\alpha u_\varepsilon(t)\|_{L^2,\varepsilon}^2$ using (\ref{eq:sysepsilon}).
First, because $A_{0\varepsilon}$ is pointwise Hermitian,
\begin{eqnarray}\label{eq:ddtuepssquare}
\nonumber\frac{d}{dt}\|\partial^\alpha u_\varepsilon(t)\|_{L^2,\varepsilon}^2&=&\frac{d}{dt}\left(A_{0\varepsilon}\partial^\alpha u_\varepsilon,\partial^\alpha u_\varepsilon\right)_{L^2}\\
\nonumber&=&\Re e\left(\frac{\partial A_{0\varepsilon}}{\partial t}\cdot\partial^\alpha u_\varepsilon,\partial^\alpha u_\varepsilon\right)_{L^2}+2\Re e\left(A_{0\varepsilon}\cdot\frac{\partial \partial^\alpha u_\varepsilon}{\partial t},\partial^\alpha u_\varepsilon\right)_{L^2}\\
\nonumber&=&\Re e\left(\frac{\partial A_{0\varepsilon}}{\partial t}\cdot\partial^\alpha u_\varepsilon,\partial^\alpha u_\varepsilon\right)_{L^2}+2\Re e\left(A_{0\varepsilon}\cdot\partial^\alpha\frac{\partial u_\varepsilon}{\partial t},\partial^\alpha u_\varepsilon\right)_{L^2}\\
\nonumber&=&\Re e\left(\frac{\partial A_{0\varepsilon}}{\partial t}\cdot\partial^\alpha u_\varepsilon,\partial^\alpha u_\varepsilon\right)_{L^2}+2\Re e\left(\partial^\alpha (A_{0\varepsilon}\cdot\frac{\partial u_\varepsilon}{\partial t}),\partial^\alpha u_\varepsilon\right)_{L^2}\\
\nonumber&&+2\Re e\left([A_{0\varepsilon},\partial^\alpha]\frac{\partial u_\varepsilon}{\partial t},\partial^\alpha u_\varepsilon\right)_{L^2}\\
\nonumber&\bui{=}{\rm(\ref{eq:sysepsilon})}&\Re e\left(\frac{\partial A_{0\varepsilon}}{\partial t}\cdot\partial^\alpha u_\varepsilon,\partial^\alpha u_\varepsilon\right)_{L^2}+2\Re e\left(\partial^\alpha J_\varepsilon L_\varepsilon J_\varepsilon u_\varepsilon,\partial^\alpha u_\varepsilon\right)_{L^2}+2\Re e\left(\partial^\alpha g_\varepsilon,\partial^\alpha u_\varepsilon\right)_{L^2}\\
&&+2\Re e\left([A_{0\varepsilon},\partial^\alpha]\frac{\partial u_\varepsilon}{\partial t},\partial^\alpha u_\varepsilon\right)_{L^2},
\end{eqnarray} 
where we have denoted $L_\varepsilon:=L(t,x,J_\varepsilon u_\varepsilon)$ and $g_\varepsilon:=J_\varepsilon g(t,x,J_\varepsilon u_\varepsilon)$.
The first term in the r.h.s. of (\ref{eq:ddtuepssquare}) can easily be estimated:
\begin{eqnarray*} 
|\Re e\left(\frac{\partial A_{0\varepsilon}}{\partial t}\cdot\partial^\alpha u_\varepsilon,\partial^\alpha u_\varepsilon\right)_{L^2}|&\leq&\|\frac{\partial A_{0\varepsilon}}{\partial t}\cdot\partial^\alpha u_\varepsilon(t)\|_{L^2}\cdot\|\partial^\alpha u_\varepsilon(t)\|_{L^2}\\
&\leq&\|\frac{\partial A_{0\varepsilon}}{\partial t}\|_{L^\infty}\cdot\|\partial^\alpha u_\varepsilon(t)\|_{L^2}^2\\
&\leq&C_1(\|J_\varepsilon u_\varepsilon(t)\|_{L^\infty},\|\frac{\partial J_\varepsilon u_\varepsilon}{\partial t}(t)\|_{L^\infty})\cdot\|\partial^\alpha u_\varepsilon(t)\|_{L^2}^2\\
&\leq&C_2(\|u_\varepsilon(t)\|_{L^\infty},\|\frac{\partial u_\varepsilon}{\partial t}(t)\|_{L^\infty})\cdot\|\partial^\alpha u_\varepsilon(t)\|_{L^2}^2\\
&\bui{\leq}{\rm(\ref{eq:sysepsilon})}&C_3(\|u_\varepsilon(t)\|_{L^\infty},\|u_\varepsilon(t)\|_{C^1})\cdot\|\partial^\alpha u_\varepsilon(t)\|_{L^2}^2\\
&\leq&C_4(\|u_\varepsilon(t)\|_{H^k})\cdot\|u_\varepsilon(t)\|_{H^k}^2,
\end{eqnarray*} 
where we have used the continuous embedding $H^k(\mathbb{T}^n)\hookrightarrow C^1(\mathbb{T}^n)$ (valid because of $k>\frac{n}{2}+1$) as well as $\|J_\varepsilon\|_{\mathcal{L}(C^l,C^l)}\leq 1$ for any $l\in\mathbb{N}$.
Let us consider the second term in the r.h.s. of (\ref{eq:ddtuepssquare}).
Since we may choose the mollifier $J_\varepsilon$ such that $J_\varepsilon^*= J_\varepsilon$ in $L^2$ (choose e.g. $\theta\in C_c^\infty(\R^n)$ with $\theta\circ(-\mathrm{Id})=\theta$), we have
\begin{eqnarray*} 
2\Re e\left(\partial^\alpha J_\varepsilon L_\varepsilon J_\varepsilon u_\varepsilon,\partial^\alpha u_\varepsilon\right)_{L^2}&=&2\Re e\left(J_\varepsilon\partial^\alpha  L_\varepsilon J_\varepsilon u_\varepsilon,\partial^\alpha u_\varepsilon\right)_{L^2}\\
&=&2\Re e\left(\partial^\alpha  L_\varepsilon J_\varepsilon u_\varepsilon,\partial^\alpha J_\varepsilon u_\varepsilon\right)_{L^2}\\
&=&2\Re e\left(  L_\varepsilon\partial^\alpha J_\varepsilon  u_\varepsilon,\partial^\alpha J_\varepsilon u_\varepsilon\right)_{L^2}+2\Re e\left([\partial^\alpha, L_\varepsilon] J_\varepsilon  u_\varepsilon,\partial^\alpha J_\varepsilon u_\varepsilon\right)_{L^2}\\
&=&\left(  (L_\varepsilon+L_\varepsilon^*) \partial^\alpha J_\varepsilon u_\varepsilon,\partial^\alpha J_\varepsilon u_\varepsilon\right)_{L^2}+2\Re e\left([\partial^\alpha, L_\varepsilon] J_\varepsilon  u_\varepsilon,\partial^\alpha J_\varepsilon u_\varepsilon\right)_{L^2},
\end{eqnarray*} 
where $L_\varepsilon^*$ is the formal adjoint of the differential operator $L_\varepsilon$.
Now, since by assumption $A_j=A_j^*$ pointwise, we have $L_\varepsilon^*=-\sum_{j=1}^n\partial_j(A_j(t,x,J_\varepsilon u_\varepsilon)\cdot)$, so that 
\[L_\varepsilon+L_\varepsilon^*=-\sum_{j=1}^n\partial_jA_j(t,x,J_\varepsilon u_\varepsilon)\]
is of zero order (this is one of the main places where symmetric hyperbolicity is used), so that 
\begin{eqnarray}\label{eq:LL*}
\nonumber|\left(  (L_\varepsilon+L_\varepsilon^*) \partial^\alpha J_\varepsilon u_\varepsilon,\partial^\alpha J_\varepsilon u_\varepsilon\right)_{L^2}|&\leq&\|(L_\varepsilon+L_\varepsilon^*) \partial^\alpha J_\varepsilon u_\varepsilon\|_{L^2}\cdot\|\partial^\alpha J_\varepsilon u_\varepsilon\|_{L^2}\\
\nonumber&\leq&\sum_{j=1}^n\|\partial_jA_j(t,x,J_\varepsilon u_\varepsilon)\|_{L^\infty}\cdot\|\partial^\alpha J_\varepsilon u_\varepsilon\|_{L^2}^2\\
\nonumber&\leq&C_5(\|J_\varepsilon u_\varepsilon(t)\|_{C^1})\cdot\|u_\varepsilon(t)\|_{H^k}^2\\
&\leq&C_6(\|u_\varepsilon(t)\|_{H^k})\cdot\|u_\varepsilon(t)\|_{H^k}^2.
\end{eqnarray} 
With 
\begin{eqnarray*} 
[\partial^\alpha, L_\varepsilon]v&=&\sum_{j=1}^n\partial^\alpha(A_j(t,x,J_\varepsilon u_\varepsilon)\partial_jv)-A_j(t,x,J_\varepsilon u_\varepsilon)\partial_j(\partial^\alpha v)\\
&=&\sum_{j=1}^n\partial^\alpha(A_j(t,x,J_\varepsilon u_\varepsilon)\partial_jv)-A_j(t,x,J_\varepsilon u_\varepsilon)\partial^\alpha(\partial_jv)\\
&=&\sum_{j=1}^n[\partial^\alpha,A_j(t,x,J_\varepsilon u_\varepsilon)]\partial_jv,
\end{eqnarray*}
we have 
\begin{eqnarray*} 
\|[\partial^\alpha, L_\varepsilon]v\|_{L^2}&\leq&\sum_{j=1}^n\|[\partial^\alpha,A_j(t,x,J_\varepsilon u_\varepsilon)]\partial_jv\|_{L^2}\\
&\bui{\leq}{\rm(\ref{eq:Moser2})}&C_7\cdot\sum_{j=1}^n\left(\|\nabla A_j(t,x,J_\varepsilon u_\varepsilon)\|_{H^{k-1}}\cdot\|\partial_jv\|_{L^\infty}+\|\nabla A_j(t,x,J_\varepsilon u_\varepsilon)\|_{L^\infty}\cdot\|\partial_jv\|_{H^{k-1}}\right)\\
&\leq&C_8\cdot\sum_{j=1}^n\left(\|A_j(t,x,J_\varepsilon u_\varepsilon)\|_{H^{k}}\cdot\|v\|_{C^1}+\|A_j(t,x,J_\varepsilon u_\varepsilon)\|_{C^1}\cdot\|v\|_{H^{k}}\right)\\
&\bui{\leq}{\rm(\ref{eq:Moser3bis})}&C_8\cdot\left(C_9(\|J_\varepsilon u_\varepsilon(t)\|_{L^\infty})\cdot(1+\|J_\varepsilon u_\varepsilon(t)\|_{H^k})\cdot\|v\|_{C^1}+C_{10}(\|J_\varepsilon u_\varepsilon(t)\|_{C^1})\cdot\|v\|_{H^{k}}\right)\\
&\leq&C_{11}(\|u_\varepsilon(t)\|_{H^k})\cdot\|v\|_{H^k},
\end{eqnarray*} 
so that 
\begin{eqnarray*} 
|2\Re e\left([\partial^\alpha, L_\varepsilon] J_\varepsilon  u_\varepsilon,\partial^\alpha J_\varepsilon u_\varepsilon\right)_{L^2}|&\leq&2\|[\partial^\alpha, L_\varepsilon] J_\varepsilon  u_\varepsilon\|_{L^2}\cdot\|\partial^\alpha J_\varepsilon u_\varepsilon\|_{L^2}\\
&\leq&C_{11}(\|u_\varepsilon(t)\|_{H^k})\cdot\|J_\varepsilon  u_\varepsilon(t)\|_{H^k}\cdot\|J_\varepsilon  u_\varepsilon(t)\|_{H^k}\\
&\leq&C_{11}(\|u_\varepsilon(t)\|_{H^k})\cdot\|u_\varepsilon(t)\|_{H^k}^2,
\end{eqnarray*} 
which gives, together with (\ref{eq:LL*}) and using $[J_\varepsilon,\partial^\alpha]=0$,
\[2|\Re e\left(\partial^\alpha J_\varepsilon L_\varepsilon J_\varepsilon u_\varepsilon,\partial^\alpha u_\varepsilon\right)_{L^2}|\leq C_{12}(\|u_\varepsilon(t)\|_{H^k}).\]
For the third term in the r.h.s. of (\ref{eq:ddtuepssquare}), we have 
\begin{eqnarray*} 
2|\Re e\left(\partial^\alpha g_\varepsilon,\partial^\alpha u_\varepsilon\right)_{L^2}|&\leq&2\|\partial^\alpha g_\varepsilon\|_{L^2}\cdot\|\partial^\alpha u_\varepsilon\|_{L^2}\\
&\bui{\leq}{\rm(\ref{eq:Moser3bis})}&C_{13}(\|u_\varepsilon(t)\|_{H^k})\cdot(1+\|J_\varepsilon  u_\varepsilon(t)\|_{H^k}\cdot\|u_\varepsilon(t)\|_{H^k})\\
&\leq&C_{14}(\|u_\varepsilon(t)\|_{H^k}).
\end{eqnarray*}
The last term in the r.h.s. of (\ref{eq:ddtuepssquare}) can be estimated as follows:
\begin{eqnarray*} 
2|\Re e\left([A_{0\varepsilon},\partial^\alpha]\frac{\partial u_\varepsilon}{\partial t},\partial^\alpha u_\varepsilon\right)_{L^2}|&\bui{=}{\rm(\ref{eq:sysepsilon})}&2|\Re e\left([A_{0\varepsilon},\partial^\alpha]A_{0\varepsilon}^{-1}(J_\varepsilon L_\varepsilon J_\varepsilon u_\varepsilon+g_\varepsilon),\partial^\alpha u_\varepsilon\right)_{L^2}|\\
&\leq&2\|[A_{0\varepsilon},\partial^\alpha]A_{0\varepsilon}^{-1}(J_\varepsilon L_\varepsilon J_\varepsilon u_\varepsilon+g_\varepsilon)\|_{L^2}\cdot\|\partial^\alpha u_\varepsilon\|_{L^2}\\
&\bui{\leq}{\rm(\ref{eq:Moser2})}&C_{15}\cdot\|u_\varepsilon\|_{H^k}\cdot\Big(\|\nabla A_{0\varepsilon}\|_{H^{k-1}}\cdot\|A_{0\varepsilon}^{-1}\cdot(J_\varepsilon L_\varepsilon J_\varepsilon u_\varepsilon+g_\varepsilon)\|_{L^\infty}\\
&&\phantom{C\cdot\|u_\varepsilon\|_{H^k}\cdot\Big(}+\|\nabla A_{0\varepsilon}\|_{L^\infty}\cdot\|A_{0\varepsilon}^{-1}\cdot(J_\varepsilon L_\varepsilon J_\varepsilon u_\varepsilon+g_\varepsilon)\|_{H^{k-1}}\Big)\\
%&\leq&C(\|u_\varepsilon\|_{H^k})\cdot\Big(C(\|u_\varepsilon\|_{C^1})+C(\|u_\varepsilon\|_{H^k})\Big)\\
&\leq&C_{16}(\|u_\varepsilon\|_{H^k}),
\end{eqnarray*}
where, in the last step, we have used the Moser estimates (\ref{eq:Moser1}) and (\ref{eq:Moser3bis}).
On the whole, $|\frac{d}{dt}\|\partial^\alpha u_\varepsilon(t)\|_{L^2,\varepsilon}^2|\leq C(\|u_\varepsilon\|_{H^k})$, so that, using the equivalence of the norms $\|\cdot\|_{L^2}$ and $\|\cdot\|_{L^2,\varepsilon}$ on some (fixed) compact set $[-T,T]\times\mathbb{T}^n$, we deduce that, setting $\|v\|_{H^k,\varepsilon}^2:=\sum_{|\alpha|\leq k}\|\partial^\alpha v\|_{L^2,\varepsilon}^2$
\[|\frac{d}{dt}\|u_\varepsilon(t)\|_{H^k,\varepsilon}^2|\leq C_{17}(\|u_\varepsilon(t)\|_{H^k,\varepsilon}).\]
By Bihari's inequality \cite{Bihari56}, we deduce that there exists a function $K$, defined and continuous on a sufficiently small interval $]-\eta,\eta[$ about $0$, such that $\|u_\varepsilon(t)\|_{H^k,\varepsilon}\leq K(t)$ for all $t\in]-\eta_\varepsilon,\eta_\varepsilon[$.
Up to making $\eta>0$ smaller, we may assume that $K(t)\leq K'<\infty$ for all $t\in]-\eta,\eta[$, so that $\|u_\varepsilon(t)\|_{H^k,\varepsilon}\leq K'$ and hence also $\|u_\varepsilon(t)\|_{H^k}\leq K$ for some $K\in]0,\infty[$.
The last statement of Claim 3 follows from the extension criterion for ODE's (valid by Steps 1 and/or 2), stating that, by  $\|u_\varepsilon(t)\|_{H^k}\leq K<\infty$, the solution $u_\varepsilon$ can be \emph{a fortiori} extended onto $]-\eta,\eta[$, QED.
\qed\\ 
\end{proof}

{\bf Step 4}: The preceding uniform estimate shows boundedness of approximate solutions in certain Sobolev spaces; use weak $*$-compactness to deduce the existence of an accumulation point. 
Then apply the interpolation inequalities (allowing compact embeddings into fractional Sobolev spaces) to deduce that the solution is actually $C^0C^1\cap C^1C^0$, that is, $C^1$ (use uniform continuity because of compactness of $[-\eta,\eta]\times\mathbb{T}^n$).\\

{\bf Claim 4}: {\em The family $(u_\varepsilon)_\varepsilon\subset C^1(]-\eta,\eta[,H^k)$ from {\rm Claim 3}, when restricted to any compact interval $I\subset]-\eta,\eta[$, admits a weak (in a particular sense) limit point $u\in C^1(I\times\mathbb{T}^n)\cap L^\infty(I,H^k)\cap{\rm Lip}(I,H^{k-1})$ which solves {\rm (\ref{eq:symmhypsys1storder})}.}\\

\begin{proof}
From Step 3 we have the existence of an $\eta>0$ and a $K\in]0,\infty[$ such that, for all sufficiently small $\varepsilon>0$, the approximate solution $u_\varepsilon$ lies in $C^1(]-\eta,\eta[,H^k)$ with $\|u_\varepsilon\|_{C^0(]-\eta,\eta[,H^k)}\leq K$.
%Moreover, there is a continuous function $K(t)$ such that $\|u_\varepsilon(t)\|_{H^k}\leq K(t)$ for all $t\in]-\eta,\eta[$.
Hence fixing an arbitrary compact interval $I\subset]-\eta,\eta[$, we have $\|u_\varepsilon\|_{C^0(I,H^k)}\leq K$, in particular the family $(u_\varepsilon)_\varepsilon$ is bounded in $C^0(I,H^k)$ and thus in $L^\infty(I,H^k)$.
Using (\ref{eq:sysepsilon}) and Moser estimates, the norm $\|\frac{\partial u_\varepsilon}{\partial t}\|_{C^0H^{k-1}}$ can be uniformly in $\varepsilon$ estimated in terms of $\|u_\varepsilon\|_{C^0(I,H^k)}$ and hence the family $(\frac{\partial u_\varepsilon}{\partial t})_\varepsilon$ is bounded in $C^0(I,H^{k-1})$, so that  $(u_\varepsilon)_\varepsilon$ is bounded in $C^1(I,H^{k-1})$ and therefore in $\mathrm{Lip}(I,H^{k-1})$.
Now $L^\infty(I,H^k)=L^1(I,H^k)'$ (topological dual), $\mathrm{Lip}(I,H^{k-1})=W^{1,\infty}(I,H^{k-1})$ by Rademacher's theorem and the latter space in turn can be identified with a closed subspace of $L^\infty(I,H^{k-1})\oplus L^\infty(I,H^{k-1})=L^1(I,H^{k-1})'\oplus L^1(I,H^{k-1})'$ via $f\mapsto (f,f')$.
Since the unit ball of the dual space of any Banach space is weakly $*$-compact, there exists a sequence $\varepsilon_p\to0$, a $u\in L^\infty(I,H^k)\cap \mathrm{Lip}(I,H^{k-1})$, such that $(u_{\varepsilon_p})_p$ converges to $u$ $*$-weakly in both spaces.
On the other hand, since $k>\frac{n}{2}+1$ and, for any $\sigma\in]0,k-\frac{n}{2}-1[$, the embedding $H^{k-\sigma}\subset C^1$ is compact, we can assume up to taking subsequences that $(u_{\varepsilon_p})_p$ converges in $C^0C^1$ to a $\overline{u}\in C^0C^1$; in fact $\overline{u}=u$ since both can be seen as sitting in the space $L^\infty(I,C^1)$ and both convergences imply the convergence in a weaker sense.
Similarly, for any $\sigma\in]0,k-\frac{n}{2}-1[$, the embedding $H^{k-1-\sigma}\subset C^0$ is compact, hence so is $C^1H^{k-1}\subset C^1C^0$, so that we may assume that $(u_{\varepsilon_p})_p$ converges in $C^1C^0$ to some $\hat{u}\in C^1C^0$ and again $\hat{u}=u$.
Since $u$ is the limit of $(u_{\varepsilon_p})_p$ in the $C^1$-topology and $J_\varepsilon\buil{\longrightarrow}{\varepsilon\to0}\mathrm{Id}$ pointwise in $C^0C^1$, we deduce that $u$ solves (\ref{eq:symmhypsys1storder}).
\qed\\
\end{proof}

{\bf Step 5}: Look at the pointwise (in $t$) $L^2$-norm of the difference between an exact $C^1$ solution to (\ref{eq:symmhypsys1storder}) and an approximate solution for any $\varepsilon>0$.
Estimate that norm on $I$ using standard estimates and Bihari's inequality.
The key point at the end is to show that $\|\mathrm{Id}-J_\varepsilon\|_{\mathcal{L}(H^1,L^2)}\leq C_1\cdot\varepsilon$ for some constant $C_1>0$.\\

{\bf Claim 5}: {\em Given any $h\in H^{k}(\mathbb{T}^n)$ and any $\varepsilon>0$, let $u_\varepsilon\in C^1(I,H^{k,2}(\mathbb{T}^n))$ solve 
\begin{equation}\label{eq:symmhypepsilon}
\left\{\begin{array}{ll}A_0(t,x,J_\varepsilon u_\varepsilon)\frac{\partial u_\varepsilon}{\partial t}&=J_\varepsilon L_\varepsilon J_\varepsilon u_\varepsilon+g_\varepsilon\textrm{ on }I\\ u_\varepsilon(0)&=h\end{array}\right. 
\end{equation}
with $u_\varepsilon$ is bounded uniformly in $\varepsilon>0$ in the $C^0H^k$-norm for all $\varepsilon$.
Let $u\in C^1$ solve {\rm(\ref{eq:symmhypsys1storder})} and consider $v_\varepsilon:=u-u_\varepsilon$.
Then there is a function $a(t):=C(\|u_\varepsilon(t)\|_{C^1},\|u(t)\|_{C^1})$ for all $t$ such that
\[\|v_\varepsilon(t)\|_{L^2}^2\leq\exp\left(\int_0^t a(s)ds\right)\cdot\left(\|\underbrace{f-h}_{v_\varepsilon(0)}\|_{L^2}^2+\int_0^t C_2(\|u_\varepsilon(s)\|_{H^k})\cdot\varepsilon\cdot e^{-\int_0^s a(\tau)d\tau}ds\right).\]
In particular, it follows from the boundedness of $(\|u_\varepsilon\|_{C^0H^k})_\varepsilon$ in Claim 3 that $u$ is unique.\\}

\begin{proof} 
We estimate $\|v_\varepsilon(t)\|_{L^2}^2$ for all $t\in I$.
First, with the notations introduced above, we write 
\begin{eqnarray}\label{eq:partialvepspartialt}
\nonumber\frac{\partial v_\varepsilon}{\partial t}&=&A_0^{-1}\frac{\partial u}{\partial t}-A_{0\varepsilon}^{-1}\frac{\partial u_\varepsilon}{\partial t}\\
\nonumber&=&A_0^{-1}L(t,x,u,\partial)u-A_{0\varepsilon}^{-1}J_\varepsilon L_\varepsilon J_\varepsilon u_\varepsilon+A_0^{-1}g(t,x,u)-A_{0\varepsilon}^{-1}g_\varepsilon\\
\nonumber&=&A_0^{-1}L(t,x,u,\partial)v_\varepsilon+(A_0^{-1}-A_{0\varepsilon}^{-1})L(t,x,u,\partial)u_\varepsilon+A_{0\varepsilon}^{-1}(L(t,x,u,\partial)u_\varepsilon-J_\varepsilon L_\varepsilon J_\varepsilon u_\varepsilon)\\
&&+(A_0^{-1}-A_{0\varepsilon}^{-1})g(t,x,u)+A_{0\varepsilon}^{-1}(g(t,x,u)-g_\varepsilon).
\end{eqnarray}
We start looking at the difference
\begin{eqnarray*} 
L(t,x,u,\partial)u_\varepsilon-J_\varepsilon L(t,x,J_\varepsilon u_\varepsilon,\partial)J_\varepsilon u_\varepsilon&=&(L(t,x,u,\partial)-L(t,x,u_\varepsilon,\partial))u_\varepsilon\\
&&+(\mathrm{Id}-J_\varepsilon)L(t,x,u_\varepsilon,\partial)u_\varepsilon+J_\varepsilon L(t,x,u_\varepsilon,\partial)(\mathrm{Id}-J_\varepsilon)u_\varepsilon\\
&&+J_\varepsilon\left(L(t,x,u_\varepsilon,\partial)-L(t,x,J_\varepsilon u_\varepsilon,\partial)\right)J_\varepsilon u_\varepsilon
\end{eqnarray*} 
and 
\begin{eqnarray*} 
g(t,x,u)-J_\varepsilon g(t,x,J_\varepsilon u_\varepsilon)&=&g(t,x,u)-g(t,x,u_\varepsilon)+g(t,x,u_\varepsilon)-g(t,x,J_\varepsilon u_\varepsilon)\\
&&+g(t,x,J_\varepsilon u_\varepsilon)-J_\varepsilon g(t,x,J_\varepsilon u_\varepsilon).
\end{eqnarray*} 
Since $A_j,g\in C^1(I\times\mathbb{T}^n)$, we may write, for all $w_1,w_2\in\mathbb{K}^N$,
\[g(t,x,w_1)-g(t,x,w_2)=\int_0^1\partial_zg(t,x,(1-s)w_2+sw_1)(w_1-w_2)ds=:G(w_1,w_2)(w_1-w_2)\]
where $\partial_z g$ denotes the derivative of $w\mapsto g(t,x,w)$ and similarly for the first-order operator
\begin{eqnarray*} 
L(t,x,w_1,\partial)-L(t,x,w_2,\partial)&=&\sum_{j=1}^n\int_0^1\partial_zA_j(t,x,(1-s)w_2+sw_1)(w_1-w_2)ds\frac{\partial}{\partial x_j}\\
&=:&M(t,x,w_1,w_2)(w_1-w_2).
\end{eqnarray*} 
In the same way, we can write 
\begin{eqnarray*} 
A_0^{-1}(t,x,u)-A_0^{-1}(t,x,J_\varepsilon u_\varepsilon)&=&A_0^{-1}(t,x,u)-A_0^{-1}(t,x,u_\varepsilon)+A_0^{-1}(t,x,u_\varepsilon)-A_0^{-1}(t,x,J_\varepsilon u_\varepsilon)\\
&=&\int_0^1d_{(t,x,(1-s)u_\varepsilon+su)}(A_0^{-1})(u-u_\varepsilon)ds\\
&&+\int_0^1d_{(t,x,(1-s)J_\varepsilon u_\varepsilon+su_\varepsilon)}(A_0^{-1})(u_\varepsilon-J_\varepsilon u_\varepsilon)ds\\
&=:&M_0(t,x,u,u_\varepsilon)(v_\varepsilon)+N_0(t,x,u_\varepsilon)(\mathrm{Id}-J_\varepsilon)(u_\varepsilon).
\end{eqnarray*} 
It is very important to notice that $G$, $M$, $M_0$ and $N_0$ depend only pointwise on $u,u_\varepsilon\ldots$, so that they can be estimated in terms of the $C^0$-norms of $u,u_\varepsilon\ldots$ only.
Now, we split the r.h.s. of (\ref{eq:partialvepspartialt}) according to their dependence on $v_\varepsilon$ and obtain 
\[\frac{\partial v_\varepsilon}{\partial t}=A_0^{-1}L(t,x,u,\partial)v_\varepsilon+A(t,x,u,u_\varepsilon,\nabla u_\varepsilon)v_\varepsilon+R_\varepsilon,\]
where
\begin{eqnarray*} 
A(t,x,u,u_\varepsilon,\nabla u_\varepsilon)&:=&M_0(t,x,u,u_\varepsilon)(v_\varepsilon)(L(t,x,u,\partial)u_\varepsilon+g(t,x,u))\\
&&+A_{0\varepsilon}^{-1}(M(t,x,u,u_\varepsilon)(v_\varepsilon)u_\varepsilon+G(u,u_\varepsilon)(v_\varepsilon)) 
\end{eqnarray*}
and
\begin{eqnarray*}
\nonumber R_\varepsilon&:=&N_0(t,x,u_\varepsilon)(\mathrm{Id}-J_\varepsilon)(u_\varepsilon)L(t,x,u,\partial)u_\varepsilon+A_{0\varepsilon}^{-1}(\mathrm{Id}-J_\varepsilon)L(t,x,u_\varepsilon,\partial)u_\varepsilon\\
\nonumber &&+A_{0\varepsilon}^{-1}J_\varepsilon L(t,x,u_\varepsilon,\partial)(\mathrm{Id}-J_\varepsilon)u_\varepsilon+A_{0\varepsilon}^{-1}J_\varepsilon M(t,x,u_\varepsilon,J_\varepsilon u_\varepsilon)(\mathrm{Id}-J_\varepsilon)(u_\varepsilon)(J_\varepsilon u_\varepsilon)\\
&&+N_0(t,x,u_\varepsilon)(\mathrm{Id}-J_\varepsilon)(u_\varepsilon)(g(t,x,u))+A_{0\varepsilon}^{-1}G(u_\varepsilon,J_\varepsilon u_\varepsilon)(\mathrm{Id}-J_\varepsilon)u_\varepsilon+A_{0\varepsilon}^{-1}(\mathrm{Id}-J_\varepsilon)g(t,x,J_\varepsilon u_\varepsilon).
\end{eqnarray*}

Next we estimate $\|R_\varepsilon(t)\|_{L^2(\mathbb{T}^n)}$ term by term.
% Using the multiplication lemma \cite[Prop. 13.3.7]{TaylorPDEIII}, stating that there exists a constant $C=C(k,n)\in\R_+^\times$ such that, for all $f,g\in L^\infty\cap H^k(\R^n)$,
% \[\|f\cdot g\|_{H^k}\leq C\cdot\left(\|f\|_{L^\infty}\|g\|_{H^k}+\|f\|_{H^k}\|g\|_{L^\infty}\right)\]
We estimate the first term as follows:
\begin{eqnarray*} 
\hspace{-1cm}\|N_0(t,x,u_\varepsilon)(\mathrm{Id}-J_\varepsilon)(u_\varepsilon)L(t,x,u_\varepsilon,\partial)u_\varepsilon\|_{L^2(\mathbb{T}^n)}(t)&\leq&\|N_0(t,x,u_\varepsilon)\|_{L^\infty}\cdot\|(\mathrm{Id}-J_\varepsilon)(u_\varepsilon)\|_{L^2}\cdot\|L(t,x,u_\varepsilon,\partial)u_\varepsilon\|_{L^\infty}\\
&\leq&C_2(\|u_\varepsilon(t)\|_{C^0})\cdot\|\mathrm{Id}-J_\varepsilon\|_{\mathcal{L}(H^{k},L^2)}\cdot\|u_\varepsilon(t)\|_{H^k}\cdot C_3(\|u_\varepsilon(t)\|_{C^1})\\
&\leq&C_4(\|u_\varepsilon(t)\|_{C^1})\cdot\|\mathrm{Id}-J_\varepsilon\|_{\mathcal{L}(H^{k},L^2)}\cdot\|u_\varepsilon(t)\|_{H^k}.
\end{eqnarray*}
For the second term
\begin{eqnarray*} 
\hspace{-2cm}\|A_{0\varepsilon}^{-1}(\mathrm{Id}-J_\varepsilon)L(t,x,u_\varepsilon,\partial)u_\varepsilon\|_{L^2(\mathbb{T}^n)}(t)&\leq&C_5(\|u_\varepsilon(t)\|_{C^0})\cdot\|\mathrm{Id}-J_\varepsilon\|_{\mathcal{L}(H^{k-1},L^2)}\cdot\|L(t,x,u_\varepsilon,\partial)u_\varepsilon\|_{H^{k-1}(\mathbb{T}^n)}(t)\\
&\bui{\leq}{\rm(\ref{eq:Moser1})}&C_5(\|u_\varepsilon(t)\|_{C^0})\cdot\|\mathrm{Id}-J_\varepsilon\|_{\mathcal{L}(H^{k-1},L^2)}\cdot\sum_{j=1}^n\|A_j(t,x,u_\varepsilon)\|_{L^\infty}\cdot\|\partial_ju_\varepsilon\|_{H^{k-1}}\\
&&\phantom{C\cdot\|\mathrm{Id}-J_\varepsilon\|_{\mathcal{L}(H^{k-1},L^2)}\cdot\sum}+\|A_j(t,x,u_\varepsilon)\|_{H^{k-1}}\cdot\|\partial_ju_\varepsilon\|_{L^\infty}\\
&\bui{\leq}{\rm(\ref{eq:Moser3bis})}&C_5(\|u_\varepsilon(t)\|_{C^0})\cdot\|\mathrm{Id}-J_\varepsilon\|_{\mathcal{L}(H^{k-1},L^2)}\cdot\sum_{j=1}^nC_6(\|u_\varepsilon(t)\|_{L^\infty})\cdot\|u_\varepsilon\|_{H^{k}}\\
&&\phantom{C\cdot\|\mathrm{Id}-J_\varepsilon\|_{\mathcal{L}(H^{k-1},L^2)}\cdot\sum}+\|u_\varepsilon(t)\|_{C^1}\cdot C_7(\|u_\varepsilon(t)\|_{L^\infty})(1+\|u_\varepsilon\|_{H^{k-1}})\\
&\leq&C_8(\|u_\varepsilon(t)\|_{C^1})\cdot \|\mathrm{Id}-J_\varepsilon\|_{\mathcal{L}(H^{k-1},L^2)}\cdot(1+\|u_\varepsilon\|_{H^{k}}).
\end{eqnarray*} 
In the same way, using also $\|J_\varepsilon\|_{\mathcal{L}(L^2,L^2)}\leq1$ and $[\partial_j,J_\varepsilon]=0$,
\begin{eqnarray*} 
\|A_{0\varepsilon}^{-1}J_\varepsilon L(t,x,u_\varepsilon,\partial)(\mathrm{Id}-J_\varepsilon)u_\varepsilon\|_{L^2}&\leq&C_9(\|u_\varepsilon(t)\|_{C^0})\cdot\|L(t,x,u_\varepsilon,\partial)(\mathrm{Id}-J_\varepsilon)u_\varepsilon\|_{L^2}\\
&\leq&C_9(\|u_\varepsilon(t)\|_{C^0})\cdot\sum_{j=1}^n\|A_j(t,x,u_\varepsilon)\|_{L^\infty}\cdot\|\partial_j(\mathrm{Id}-J_\varepsilon)u_\varepsilon\|_{L^2}\\
&\leq&C_{10}(\|u_\varepsilon(t)\|_{C^0})\cdot\|\mathrm{Id}-J_\varepsilon\|_{\mathcal{L}(H^{k-1},L^2)}\cdot\|u_\varepsilon\|_{H^k}
\end{eqnarray*} 
and, as $\|M(t,x,u_\varepsilon,J_\varepsilon u_\varepsilon)(u_\varepsilon-J_\varepsilon u_\varepsilon)\|_{L^2}\leq C_{11}(\|u_\varepsilon\|_{L^\infty})\cdot\|u_\varepsilon-J_\varepsilon u_\varepsilon\|_{L^2}$, we obtain
\begin{eqnarray*} 
\|A_{0\varepsilon}^{-1}J_\varepsilon M(t,x,u_\varepsilon,J_\varepsilon u_\varepsilon)(u_\varepsilon-J_\varepsilon u_\varepsilon)J_\varepsilon u_\varepsilon\|_{L^2}&\leq&C_{12}(\|u_\varepsilon(t)\|_{C^0})\cdot\|M(t,x,u_\varepsilon,J_\varepsilon u_\varepsilon)(u_\varepsilon-J_\varepsilon u_\varepsilon)\cdot J_\varepsilon u_\varepsilon\|_{L^2}\\
&\leq&C_{13}(\|u_\varepsilon(t)\|_{C^0})\cdot\|(\mathrm{Id}-J_\varepsilon)u_\varepsilon\|_{L^2}\\
&\leq&C_{13}(\|u_\varepsilon(t)\|_{C^1})\cdot\|\mathrm{Id}-J_\varepsilon\|_{\mathcal{L}(H^{k},L^2)}\cdot\|u_\varepsilon\|_{H^k}.
\end{eqnarray*} 
As before, estimating $N_0$ and $g$, we have
\begin{eqnarray*} 
\|N_0(t,x,u_\varepsilon)(\mathrm{Id}-J_\varepsilon)(u_\varepsilon)(g(t,x,u))\|_{L^2}&\leq&\|N_0(t,x,u_\varepsilon)\|_{L^\infty}\cdot\|(\mathrm{Id}-J_\varepsilon)(u_\varepsilon)\|_{L^2}\cdot\|g(t,x,u)\|_{L^\infty}\\
&\leq&C_{14}(\|u_\varepsilon(t)\|_{C^0})\cdot\|\mathrm{Id}-J_\varepsilon\|_{\mathcal{L}(H^{k},L^2)}\cdot\|u_\varepsilon\|_{H^k}.
\end{eqnarray*} 
For the last two terms, we obtain 
\begin{eqnarray*} 
\|A_{0\varepsilon}^{-1}G(u_\varepsilon,J_\varepsilon u_\varepsilon)(u_\varepsilon-J_\varepsilon u_\varepsilon)\|_{L^2}&\leq&C_{15}(\|u_\varepsilon(t)\|_{L^\infty})\cdot\|u_\varepsilon-J_\varepsilon u_\varepsilon\|_{L^2}\\
&\leq&C_{15}(\|u_\varepsilon(t)\|_{C^1})\cdot\|\mathrm{Id}-J_\varepsilon\|_{\mathcal{L}(H^{k},L^2)}\cdot\|u_\varepsilon\|_{H^k}
\end{eqnarray*} 
and
\begin{eqnarray*} 
\|A_{0\varepsilon}^{-1}(\mathrm{Id}-J_\varepsilon)g(t,x,J_\varepsilon u_\varepsilon)\|_{L^2}&\leq&C_{16}(\|u_\varepsilon(t)\|_{L^\infty})\cdot\|\mathrm{Id}-J_\varepsilon\|_{\mathcal{L}(H^{k},L^2)}\cdot\|g(t,x,J_\varepsilon u_\varepsilon)\|_{H^k}\\
&\bui{\leq}{\rm(\ref{eq:Moser3bis})}&C_{17}(\|u_\varepsilon(t)\|_{L^\infty})\cdot\|\mathrm{Id}-J_\varepsilon\|_{\mathcal{L}(H^{k},L^2)}\cdot (1+\|u_\varepsilon(t)\|_{H^k}).
\end{eqnarray*} 
Note that $\|\mathrm{Id}-J_\varepsilon\|_{\mathcal{L}(H^{k},L^2)}\leq\|\mathrm{Id}-J_\varepsilon\|_{\mathcal{L}(H^{k-1},L^2)}$.
% in particular 
% \begin{eqnarray*} 
% \|A_0^{-1}(t,x,u)-A_0^{-1}(t,x,J_\varepsilon u_\varepsilon)\|_{L^2}&\leq&C(\|u(t)\|_{C^0},\|u_\varepsilon(t)\|_{C^0})\cdot\|v_\varepsilon(t)\|_{L^2}\\
% &&+C(\|u_\varepsilon(t)\|_{C^0})\cdot\|\mathrm{Id}-J_\varepsilon\|_{\mathcal{L}(H^k,L^2)}\cdot\|u_\varepsilon(t)\|_{H^k}\\
% &\leq&C(\|u(t)\|_{C^0},\|u_\varepsilon(t)\|_{C^0})\cdot(\|v_\varepsilon(t)\|_{L^2}+1)
% \end{eqnarray*} 
On the whole, we obtain
\[\|R_\varepsilon(t)\|_{L^2}\leq C_{18}(\|u_\varepsilon(t)\|_{C^1})\cdot(1+\|u_\varepsilon(t)\|_{H^k})\cdot\|\mathrm{Id}-J_\varepsilon\|_{\mathcal{L}(H^{k-1},L^2)}.\]
We deduce that
\begin{eqnarray*} 
\frac{d}{dt}\|v_\varepsilon(t)\|_{L^2,\varepsilon}^2&=&\left(\frac{\partial A_{0\varepsilon}}{\partial t}\cdot v_\varepsilon,v_\varepsilon\right)_{L^2}+2\Re e\left(A_{0\varepsilon}\frac{\partial v_\varepsilon}{\partial t},v_\varepsilon\right)_{L^2}\\
&=&\left(\frac{\partial A_{0\varepsilon}}{\partial t}\cdot v_\varepsilon,v_\varepsilon\right)_{L^2}\\
&&+2\Re e\left\{(L(t,x,u,\partial)v_\varepsilon,v_\varepsilon)_{L^2}+(A(t,x,u,u_\varepsilon,\nabla u_\varepsilon)v_\varepsilon,v_\varepsilon)_{L^2,\varepsilon}+(R_\varepsilon,v_\varepsilon)_{L^2,\varepsilon}\right\},
\end{eqnarray*} 
with $|\left(\frac{\partial A_{0\varepsilon}}{\partial t}\cdot v_\varepsilon,v_\varepsilon\right)_{L^2}|\leq C_{19}(\|u_\varepsilon(t)\|_{C^1})\cdot\|v_\varepsilon(t)\|_{L^2}^2$ and
\begin{eqnarray*} 
\|A(t,x,u,u_\varepsilon,\nabla u_\varepsilon)v_\varepsilon\|_{L^2,\varepsilon}&\leq&\|M_0(t,x,u,u_\varepsilon)(v_\varepsilon)(L(t,x,u,\partial)u_\varepsilon+g(t,x,u))\|_{L^2,\varepsilon}\\
&&+\|A_{0\varepsilon}^{-1}(M(t,x,u,u_\varepsilon)(v_\varepsilon)u_\varepsilon+G(u,u_\varepsilon)(v_\varepsilon))\|_{L^2,\varepsilon}\\
&\leq&C_{20}(\|u(t)\|_{C^0},\|u_\varepsilon(t)\|_{C^1})\cdot\|v_\varepsilon(t)\|_{L^2}\\
&&+C_{21}(\|u_\varepsilon(t)\|_{C^0})\cdot(\|M(t,x,u,u_\varepsilon)(v_\varepsilon)u_\varepsilon\|_{L^2}+\|G(u,u_\varepsilon)v_\varepsilon(t)\|_{L^2})\\
&\leq&C_{20}(\|u(t)\|_{C^0},\|u_\varepsilon(t)\|_{C^1})\cdot\|v_\varepsilon(t)\|_{L^2}\\
&&+C_{22}(\|u_\varepsilon(t)\|_{C^0},\|u(t)\|_{C^0})\cdot\|u_\varepsilon(t)\|_{C^1}\cdot\|v_\varepsilon\|_{L^2}+C_{23}(\|u_\varepsilon(t)\|_{C^0},\|u(t)\|_{C^0})\cdot\|v_\varepsilon(t)\|_{L^2}\\
&\leq&C_{24}(\|u_\varepsilon(t)\|_{C^1},\|u(t)\|_{C^0})\cdot\|v_\varepsilon(t)\|_{L^2,\varepsilon}
\end{eqnarray*} 
as well as 
\begin{eqnarray*} 
\Re e\left(L(t,x,u,\partial)v_\varepsilon,v_\varepsilon\right)_{L^2}&=&\left((L+L^*)(t,x,u,\partial)v_\varepsilon,v_\varepsilon\right)_{L^2}\\
&=&-\sum_{j=1}^n\left(\partial_jA_j(t,x,u)\cdot v_\varepsilon,v_\varepsilon\right)_{L^2}
\end{eqnarray*} 
because $A_j^*=A_j$, so that
\[|\Re e\left(L(t,x,u,\partial)v_\varepsilon,v_\varepsilon\right)_{L^2}\leq C_{25}(\|u_\varepsilon(t)\|_{C^1})\cdot\|v_\varepsilon\|_{L^2,\varepsilon}^2.\]
Noticing that $2\Re e\left(R_\varepsilon,v_\varepsilon\right)_{L^2,\varepsilon}\leq\|R_\varepsilon\|_{L^2,\varepsilon}^2+\|v_\varepsilon\|_{L^2,\varepsilon}^2$, we obtain
\[\frac{d}{dt}\|v_\varepsilon(t)\|_{L^2,\varepsilon}^2\leq C_{26}(\|u_\varepsilon(t)\|_{C^1},\|u(t)\|_{C^1})\cdot\|v_\varepsilon\|_{L^2,\varepsilon}^2+\|R_\varepsilon(t)\|_{L^2,\varepsilon}^2.\]
Then by Gr\"onwall's lemma, for every $t\in I\cap\R_+$,
\[\|v_\varepsilon(t)\|_{L^2,\varepsilon}^2\leq\exp\left(\int_0^t a(s)ds\right)\cdot\left(\|\underbrace{f-h}_{v_\varepsilon(0)}\|_{L^2,\varepsilon}^2+\int_0^t\|R_\varepsilon(s)\|_{L^2,\varepsilon}^2\cdot e^{-\int_0^s a(\tau)d\tau}ds\right).\]
% We deduce that, for every $t\in I\cap\R_+$,
% \[\|v_\varepsilon(t)\|_{L^2}^2\leq \exp\left(\int_0^t a(s)ds\right)\cdot\left(\|u_0-h\|_{L^2}^2+\int_0^t\|R_\varepsilon(s)\|_{L^2}^2 e^{-\int_0^s a(\tau)d\tau}ds\right).\]
Note that, since $\|u_\varepsilon(t)\|_{C^1}\leq C_{27}\cdot\|u_\varepsilon(t)\|_{H^k}\leq C_{28}$ uniformly in $\varepsilon$ and $t\in I$, we may choose $a$ to be constant.
For the same reason (and by the estimate above), $\|R_\varepsilon(t)\|_{L^2,\varepsilon}^2\leq C_{29}\cdot\|\mathrm{Id}-J_\varepsilon\|_{\mathcal{L}(H^{k-1},L^2)}^2$.
This combined with the choice $h:=f$ and the equivalence of the norms $\|\cdot\|_{L^2}$ and $\|\cdot\|_{L^2,\varepsilon}$ yields
\begin{eqnarray*} 
\|v_\varepsilon(t)\|_{L^2}^2&\leq& e^{at}\cdot\int_0^t\|R_\varepsilon(s)\|_{L^2,\varepsilon}^2 e^{-as}ds\\
&\leq& C_{29} \cdot a^{-1}\cdot(e^{at}-1)\cdot \|\mathrm{Id}-J_\varepsilon\|_{\mathcal{L}(H^{1},L^2)}^2.
\end{eqnarray*}
It follows from the Proposition \ref{p:IminusJeps} that $\|v_\varepsilon(t)\|_{L^2}^2\leq C_{30}\cdot(e^{at}-1)\cdot\varepsilon$ (recall that $k-1>\frac{n}{2}\geq\frac{1}{2}$).
This implies on the one hand that any solution to (\ref{eq:symmhypsys1storder}) -- with given initial condition $f$ -- is the pointwise (in $t$) limit when $\varepsilon\to0$ of the uniquely determined family $(u_\varepsilon)_\varepsilon$, so that any two such solutions must coincide on their common interval of definition.
On the other hand, this inequality gives the $C^0L^2$-rate of convergence for $(u_\varepsilon)_\varepsilon$ to $u$.
\qed\\ 
\end{proof}

{\bf Step 6}: By what seems to be a well-known result from functional analysis (see e.g. \cite[Lemma 4.1]{Abels09}), the fact that the solution to (\ref{eq:symmhypsys1storder}) belongs to certain Sobolev spaces implies its continuity $I\to H^k$ where $H^k$ is endowed with the weak topology.
To show the strong continuity of the solution, it suffices to show the continuity of its (pointwise) $H^k$-norm.
Estimate the $t$-derivative of that norm by inserting again a $J_\varepsilon$, using standard estimates, Gr\"onwall and making $\varepsilon\to0$ to show that the norm of the solution is actually Lipschitz.\\

{\bf Claim 6}: {\em The solution $u$ from Claim 4 actually lies not only in $L^\infty(I,H^k)$ as proven in Claim 4 but also in $C^0(I,H^k)$.}\\

\begin{proof} So the continuity in the weak sense follows from \cite[Lemma 4.1]{Abels09} applied to $Y=H^k$ and $X=H^{k-1}$ (note that $Y\subset X$ densely and so does $X'\subset Y'$).
%use again Moser and probably the fact that $A_j=A_j^*$ to control $\frac{d}{dt}\|J_\varepsilon u(t)\|_{H^k}^2$ in terms of a continuous function of $\|u(t)\|_{C^1}$ and $\|u(t)\|_{H^k}^2$ itself.
To show the strong continuity, it suffices to show that $t\mapsto\|u(t)\|_{H^k}$ is continuous.
Note here that one cannot directly estimate $\frac{d}{dt}\|u(t)\|_{H^k}^2$ as before since the differential operator $L$ does not preserve $H^k$.
As in the proof of \cite[Prop. 16.1.4]{Taylor3}, we avoid this difficulty by inserting a $J_\varepsilon$ before $u$.
Setting $\left(\cdot\,,\cdot\right)_{L^2,0}:=\left(A_0\cdot\,,\cdot\right)_{L^2}$, we pick any multiindex $\alpha$ with $|\alpha|\leq k$.
Recalling that $u\in\mathrm{Lip}(I,H^{k-1})$, we have $J_\varepsilon u\in\mathrm{Lip}(I,H^k)$, in particular, the function $t\mapsto \|\partial^\alpha J_\varepsilon u(t)\|_{L^2,0}^2$ is differentiable almost everywhere.
% We derive the function $t\mapsto \|\partial^\alpha J_\varepsilon u(t)\|_{L^2,0}^2$ first in a purely formal way and make use of the following\\
% {\bf Claim:} {\it Let $T$ be any scalar-valued distribution on an open interval such that $T'=h$ for some $L^\infty$ function $h$. 
% Then $T$ is actually an almost everywhere differentiable function.}\\
% {\it Proof of Claim}: First, given any $x_0$ in the interval, the function $g(x):=\int_{x_0}^x h(t)dt$ is well-defined and almost everywhere differentiable {\red\bf (Is it sure? Check)} with almost-everywhere-defined derivative $h$.
% Therefore $(T-g)'=0$ as a distribution.
% Now the only distributions $F$ with vanishing derivative on an interval are the constant functions.
% Namely if $F'=0$, then $F(\varphi')=0$ for all $\varphi\in\mathcal{D}(I)$, where $\mathcal{D}(I)$ is the space of all compactly supported smooth functions on $I$.
% But the space $\left\{\varphi'\,|\,\varphi\in\mathcal{D}(I)\right\}$ coincides with $\left\{\varphi\in\mathcal{D}(I)\,|\,\int_I\varphi=0\right\}$.
% %{\red\bf (Easy exercise, I can write down the details if needed. N.)}
% Pick any $\theta\in\mathcal{D}(I)$ with $\int_I\theta(t)dt=1$.
% Then for any $\varphi\in\mathcal{D}(I)$, one can write $\varphi=\varphi-(\int_I\varphi)\theta+(\int_I\varphi)\theta$, where $\int_I\left\{\varphi-(\int_I\varphi)\theta\right\}=0$, so that 
% \begin{eqnarray*} 
% F(\varphi)&=&F(\varphi-(\int_I\varphi)\theta)+F((\int_I\varphi)\theta)\\
% &=&0+(\int_I\varphi)F(\theta),
% \end{eqnarray*}
% that is, $F=F(\theta)\in\mathbb{R}$.
% This proves the claim. 
% \hfill$\surd$\\
We start computing the derivative of $t\mapsto \|\partial^\alpha J_\varepsilon u(t)\|_{L^2,0}^2$:
\begin{eqnarray*} 
\hspace{-2cm}\frac{d}{dt}\|\partial^\alpha J_\varepsilon u(t)\|_{L^2,0}^2&=&\left(\frac{\partial A_{0}}{\partial t}\cdot\partial^\alpha J_\varepsilon u ,\partial^\alpha J_\varepsilon u\right)_{L^2}+2\Re e\left(\frac{\partial}{\partial t}\partial^\alpha J_\varepsilon u,\partial^\alpha J_\varepsilon u\right)_{L^2,0}\\
&=&\left(\frac{\partial A_{0}}{\partial t}\cdot\partial^\alpha J_\varepsilon u ,\partial^\alpha J_\varepsilon u\right)_{L^2}+2\Re e\left(\partial^\alpha J_\varepsilon \frac{\partial u}{\partial t},\partial^\alpha J_\varepsilon u\right)_{L^2,0}\\
&=&\left(\frac{\partial A_{0}}{\partial t}\cdot\partial^\alpha J_\varepsilon u ,\partial^\alpha J_\varepsilon u\right)_{L^2}+2\Re e\left(\partial^\alpha J_\varepsilon A_{0}\frac{\partial u}{\partial t},\partial^\alpha J_\varepsilon u\right)_{L^2}\\
&&+2\Re e\left([A_{0},J_\varepsilon]\partial^\alpha\frac{\partial u}{\partial t},\partial^\alpha J_\varepsilon u\right)_{L^2}+2\Re e\left(J_\varepsilon [A_{0},\partial^\alpha]\frac{\partial u}{\partial t},\partial^\alpha J_\varepsilon u\right)_{L^2}\\
&\bui{=}{\rm(\ref{eq:symmhypsys1storder})}&\left(\frac{\partial A_{0}}{\partial t}\cdot\partial^\alpha J_\varepsilon u ,\partial^\alpha J_\varepsilon u\right)_{L^2}+2\Re e\left([A_{0},J_\varepsilon]\partial^\alpha\frac{\partial u}{\partial t},\partial^\alpha J_\varepsilon u\right)_{L^2}+2\Re e\left(J_\varepsilon [A_{0},\partial^\alpha]\frac{\partial u}{\partial t},\partial^\alpha J_\varepsilon u\right)_{L^2}\\
&&+2\Re e\left(\partial^\alpha J_\varepsilon L(t,x,u,\partial)u,\partial^\alpha J_\varepsilon u\right)_{L^2}+2\Re e\left(\partial^\alpha J_\varepsilon g(t,x,u),\partial^\alpha J_\varepsilon u\right)_{L^2}.
\end{eqnarray*}
%{\red\bf As in the proof of \cite[Prop. 16.1.4]{Taylor3}, we implicitly assume that $t\mapsto\|\partial^\alpha J_\varepsilon u(t)\|_{L^2}$ is differentiable $I\to \R$, or at least almost everywhere; SOLVE THE EXERCISE: ANY DISTRIBUTION WITH VANISHING DERIVATIVE IS LOCALLY CONSTANT.}\\
We begin with estimating the last term.
First, if $|\alpha|\geq1$, we have 
\begin{eqnarray*} 
|2\Re e\left(\partial^\alpha J_\varepsilon g(t,x,u),\partial^\alpha J_\varepsilon u\right)_{L^2}|&\leq&2\|\partial^\alpha J_\varepsilon g(t,x,u)\|_{L^2}\cdot\|\partial^\alpha J_\varepsilon u\|_{L^2}\\
&\leq&2\|\partial^\alpha g(t,x,u)\|_{L^2}\cdot\|\partial^\alpha u\|_{L^2}\\
&\bui{\leq}{\rm(\ref{eq:Moser3})}&C_1(\|u(t)\|_{C^0})\cdot\|\partial^\alpha u\|_{L^2}^2\\
&\leq&C_2(\|u(t)\|_{C^1})\cdot\|u(t)\|_{H^k}^2.
\end{eqnarray*} 
% There is however a problem with the case $\alpha=0$ since we have not assumed $g(t,x,u)=g(u)$ as well as $g(0)=0$.
% If we do, then again (\ref{eq:Moser3}) implies 
% \[|2\Re e\left(J_\varepsilon g(u),J_\varepsilon u\right)_{L^2}|\leq C(\|u(t)\|_{C^1})\cdot\|u(t)\|_{L^2}^2\leq C(\|u(t)\|_{C^1})\cdot\|u(t)\|_{H^k}^2,\]
% otherwise -- without any further assumption on $g$ -- 
For $\alpha=0$, we can only apply (\ref{eq:Moser3bis}) and obtain
\begin{eqnarray*} 
|2\Re e\left(J_\varepsilon g(t,x,u),J_\varepsilon u\right)_{L^2}|&\leq& C_3(\|u(t)\|_{C^1})\cdot\left(1+\|u(t)\|_{L^2}\right)\cdot\|u(t)\|_{L^2}\\
&\leq&C_4(\|u(t)\|_{C^1})\cdot\left(1+\|u(t)\|_{L^2}^2\right),
\end{eqnarray*}
which actually suffices for the proof of Step 6 (as well as for the extension criterion in Step 7).
Next we decompose the last but one term as follows:
\begin{eqnarray*} 
2\Re e\left(\partial^\alpha J_\varepsilon L(t,x,u,\partial)u,\partial^\alpha J_\varepsilon u\right)_{L^2}&=&2\Re e\left(J_\varepsilon\partial^\alpha L(t,x,u,\partial)u,\partial^\alpha J_\varepsilon u\right)_{L^2}\\
&=&2\Re e\left(J_\varepsilon L\partial^\alpha u,\partial^\alpha J_\varepsilon u\right)_{L^2}+2\Re e\left(J_\varepsilon [\partial^\alpha,L] u,\partial^\alpha J_\varepsilon u\right)_{L^2}.
\end{eqnarray*} 
The second term on the r.h.s. can be easily estimated with the help of Moser estimates:
\begin{eqnarray*} 
|2\Re e\left(J_\varepsilon [\partial^\alpha,L] u,\partial^\alpha J_\varepsilon u\right)_{L^2}|&\leq&2\|J_\varepsilon [\partial^\alpha,L] u\|_{L^2}\cdot\|\partial^\alpha J_\varepsilon u\|_{L^2}\\
&\leq&2\|[\partial^\alpha,L] u\|_{L^2}\cdot\|\partial^\alpha u\|_{L^2}\\
&\bui{\leq}{\rm(\ref{eq:Moser2})}&C_5\cdot\left(\sum_{j=1}^n\|\nabla A_j\|_{L^\infty}\cdot\|u(t)\|_{H^k}+\|\nabla A_j\|_{H^{k-1}}\cdot\|\nabla u\|_{L^\infty}\right)\|u(t)\|_{H^k}\\
&\leq&C_5\cdot\left(C_6(\|u(t)\|_{C^1})\|u(t)\|_{H^k}+C_7(\|u(t)\|_{C^1})\|u(t)\|_{H^k}\|u(t)\|_{C^1}\right)\|u(t)\|_{H^k}\\
&\leq&C_8(\|u(t)\|_{C^1})\cdot\|u(t)\|_{H^k}^2.
\end{eqnarray*} 
Using symmetric hyperbolicity, we may estimate the term 
\begin{eqnarray*} 
2\Re e\left(J_\varepsilon L\partial^\alpha u,\partial^\alpha J_\varepsilon u\right)_{L^2}&=&2\Re e\left(LJ_\varepsilon\partial^\alpha u,\partial^\alpha J_\varepsilon u\right)_{L^2}+2\Re e\left([J_\varepsilon,L]\partial^\alpha u,\partial^\alpha J_\varepsilon u\right)_{L^2}\\
&=&\Re e\left((L+L^*)\partial^\alpha J_\varepsilon u,\partial^\alpha J_\varepsilon u\right)_{L^2}+2\Re e\left([J_\varepsilon,L]\partial^\alpha u,\partial^\alpha J_\varepsilon u\right)_{L^2}\\
&=&-\sum_{j=1}^n\Re e\left(\frac{\partial A_j}{\partial x_j}\cdot\partial^\alpha J_\varepsilon u,\partial^\alpha J_\varepsilon u\right)_{L^2}+2\Re e\left([J_\varepsilon,L]\partial^\alpha u,\partial^\alpha J_\varepsilon u\right)_{L^2},
\end{eqnarray*} 
where
\begin{eqnarray*} 
|\sum_{j=1}^n\Re e\left(\frac{\partial A_j}{\partial x_j}\cdot\partial^\alpha J_\varepsilon u,\partial^\alpha J_\varepsilon u\right)_{L^2}|&\leq&\sum_{j=1}^n\|\frac{\partial A_j}{\partial x_j}\cdot\partial^\alpha J_\varepsilon u\|_{L^2}\|\cdot\|\partial^\alpha J_\varepsilon u\|_{L^2}\\
&\leq&\sum_{j=1}^n\|\frac{\partial A_j}{\partial x_j}\|_{L^\infty}\cdot\|\partial^\alpha J_\varepsilon u\|_{L^2}^2\\
&\leq&C_9(\|u(t)\|_{C^1})\cdot\|u(t)\|_{H^k}^2
\end{eqnarray*} 
and, with $[J_\varepsilon,L]v=\sum_{j=1}^n[J_\varepsilon,A_j]\frac{\partial v}{\partial x_j}$,
\begin{eqnarray*} 
|2\Re e\left([J_\varepsilon,L]\partial^\alpha u,\partial^\alpha J_\varepsilon u\right)_{L^2}|&\leq&2\sum_{j=1}^n\|[J_\varepsilon,A_j]\frac{\partial \partial^\alpha u}{\partial x_j}\|_{L^2}\cdot\|\partial^\alpha J_\varepsilon u\|_{L^2}\\
&\bui{\leq}{{\rm (Lemma \ref{l:mollifiercommutator}}.iii){\rm)}}&C_{10}\cdot\sum_{j=1}^n\|A_j\|_{C^1}\cdot\|\partial^\alpha u\|_{L^2}^2\\
&\leq&C_{11}(\|u(t)\|_{C^1})\cdot\|u(t)\|_{H^k}^2.
\end{eqnarray*} 
We also have $|\left(\frac{\partial A_{0}}{\partial t}\cdot\partial^\alpha J_\varepsilon u ,\partial^\alpha J_\varepsilon u\right)_{L^2}|\leq C_{12}(\|u(t)\|_{C^1})\cdot\|u(t)\|_{H^k}^2$,
\begin{eqnarray*} 
\hspace{-1cm}|2\Re e\left(J_\varepsilon [A_{0},\partial^\alpha]\frac{\partial u}{\partial t},\partial^\alpha J_\varepsilon u\right)_{L^2}|&\leq&2\|[A_{0},\partial^\alpha]\frac{\partial u}{\partial t}\|_{L^2}\cdot\|u(t)\|_{H^k}\\
&\bui{\leq}{\rm(\ref{eq:Moser2})}&C_{13}\cdot\left(\|\nabla A_0\|_{H^{k-1}}\cdot\|\frac{\partial u}{\partial t}\|_{L^\infty}+\|\nabla A_0\|_{L^\infty}\cdot\|\frac{\partial u}{\partial t}\|_{H^{k-1}}\right)\cdot\|u(t)\|_{H^k}\\
&\bui{\leq}{\rm(\ref{eq:Moser3bis})}&\left(C_{14}(\|u(t)\|_{C^0})\cdot\|u(t)\|_{H^k}\cdot C_{15}(\|u(t)\|_{C^1})+C_{16}(\|u(t)\|_{C^1})\cdot(1+\|u(t)\|_{H^k})\right)\cdot\|u(t)\|_{H^k}\\
&\leq&C_{17}(\|u(t)\|_{C^1})\cdot\|u(t)\|_{H^k}\cdot(1+\|u(t)\|_{H^k})
\end{eqnarray*} 
and 
\begin{eqnarray*} 
|2\Re e\left([A_{0},J_\varepsilon]\partial^\alpha\frac{\partial u}{\partial t},\partial^\alpha J_\varepsilon u\right)_{L^2}|&\leq&2\|[A_{0},J_\varepsilon]\partial^\alpha\frac{\partial u}{\partial t}\|_{L^2}\cdot\|u(t)\|_{H^k}\\
&\bui{\leq}{\rm(\ref{l:mollifiercommutator})}&C_{18}(\|u(t)\|_{C^1})\cdot\|u(t)\|_{H^k}\cdot(1+\|u(t)\|_{H^k}).
%&\leq&C\cdot\|A_0\|_{C^1}\cdot\|\partial^\alpha\frac{\partial u}{\partial t}\|_{L^2}\cdot\|u(t)\|_{H^k}\\
\end{eqnarray*} 
Bringing everything together, we deduce that, setting $\|v\|_{H^k,0}^2:=\sum_{|\alpha|\leq k}\|\partial^\alpha v\|_{L^2,0}^2$,
%and using the equivalence of $\|\cdot\|_{L^2}$ and $\|\cdot\|_{L^2,0}$
% \[\frac{d}{dt}\|J_\varepsilon u(t)\|_{H^k}^2 \leq C(\|u(t)\|_{C^1})\cdot\|u(t)\|_{H^k}\cdot(1+\|u(t)\|_{H^k}),\]
% from which
\begin{equation}\label{eq:ddtJepsilonu} 
\frac{d}{dt}\|J_\varepsilon u(t)\|_{H^k,0}^2 \leq C_{19}(\|u(t)\|_{C^1})\cdot\left(1+\|u(t)\|_{H^k}^2\right).
\end{equation}
That inequality \emph{does not depend} on $\varepsilon>0$.
%{\bf\red (Check; one must probably use the fact that $u$ is the $*$-weak limit of the family $u_\varepsilon$, which is bounded in $C^0H^k$ by Step 3.)}
Since by construction of $u$ we have the existence of a constant $C_{20}$ such that $\|u(t)\|_{H^k}\leq C_{20}$ for all $t\in I$ (because $u\in L^\infty H^k$) and since $\|\cdot\|_{L^2}$ and $\|\cdot\|_{L^2,0}$ are equivalent, we deduce that $J_\varepsilon u:I\to H^k$ is $C_{21}$-Lipschitz continuous for a constant $C_{21}>0$ independent of $\varepsilon$.
Since for all $t\in I$ one has $(J_\varepsilon u)(t)\buil{\longrightarrow}{\varepsilon\to 0}u(t)$ (in the strong $H^k$-topology) and because the pointwise limit of a family of $C_{21}$-Lipschitz continuous family is again $C_{21}$-Lipschitz continuous, we obtain that $u:I\to H^k$ is $C_{21}$-Lipschitz continuous, in particular $u\in C^0H^k$.
\qed\\ 
\end{proof}

{\bf Step 7}: Use in fact the preceding estimate of the $t$-derivative of $\|J_\varepsilon u(t)\|_{H^k}$ to deduce, using Gr\"onwall and after letting $\varepsilon\to0$, that $\|u(t)\|_{H^k}$ can be controlled in terms of a continuous function of $\|u(t)\|_{C^1}$.
Conclude the proof of Theorem \ref{t:exsymmhypsys1storder}.\\

{\bf Claim 7}: {\em The solution $u\in C^0H^k$ constructed above exists as long as $\|u(t)\|_{C^1}$ remains bounded: if, for a given $T\in(0,\infty)$, there is a constant $C$ such that $\|u(t)\|_{C^1}\leq C$ for all $t\in[0,T[$, then there exists a $\delta>0$ such that the solution $u$ can be extended to a solution in $C^0([0,T+\delta],H^k)$.}\\

\begin{proof} 
%PROOF SHOULD BE CORRECTED: ESTIMATE $H^k,0$-norm in terms of $H^k$-norm\footnote{I'll do it later (N)}\\
Since by assumption $\|u(t)\|_{C^1}\leq C<\infty$ for all $t\in[0,T[$, there exists a constant $C'$ such that $C'^{-1}\|u(t)\|_{H^k}\leq\|u(t)\|_{H^k,0}\leq C'\|u(t)\|_{H^k}$ for all $t\in[0,T[$ and inequality (\ref{eq:ddtJepsilonu}) yields
\[\frac{d}{dt}\|J_\varepsilon u(t)\|_{H^k,0}^2 \leq C_1\cdot(1+\|u(t)\|_{H^k,0}^2),\]
which can be rewritten in integral form: for every $\tau>0$,
%\footnote{Choosing $\tau<0$ would lead to the same inequality $\frac{d}{dt}\|u(t)\|_{H^k,0}^2 \leq C\cdot\|u(t)\|_{H^k,0}^2$.}
\[
\frac{\|(J_\varepsilon u)(t+\tau)\|_{H^k,0}^2-\|(J_\varepsilon u)(t)\|_{H^k,0}^2}{\tau}=\frac{1}{\tau}\int_0^\tau\frac{d}{ds}\|J_\varepsilon u(s)\|_{H^k,0}^2ds\leq \frac{C_1}{\tau}\cdot\int_0^\tau1+\|u(s)\|_{H^k,0}^2ds.
\]
Using the pointwise convergence $J_\varepsilon\buil{\longrightarrow}{\varepsilon\to0}\mathrm{Id}$ in $H^k$ (and $H_{,0}^k$) and letting then $\tau\to 0^+$ lead to
\[\frac{d}{dt}\|u(t)\|_{H^k,0}^2 \leq C_1\cdot(1+\|u(t)\|_{H^k,0}^2)\]
and therefore $\|u(t)\|_{H^k,0}^2\leq(1+\|u(0)\|_{H^k,0}^2)\cdot e^{C_1t}-1$ for all $t\in[0,T[$, in particular there is a constant $K>0$ with $\|u(t)\|_{H^k}\leq K<\infty$ for all $t\in[0,T[$.
The latter inequality implies that $u$ can be extended beyond $T$, namely as follows.
Consider a small interval of the form $]T-\hat{\eta},T+\hat{\eta}[$ for some $\hat{\eta}>0$.
Because $A_j$ and $g$ are continuous, satisfy the ``strong'' local Lipschitz condition and because the time of existence for solutions to ODE's depends \emph{continuously} on the norm of the initial condition (see e.g. proof of \cite[Theorem 6.2.1]{Pazy}), up to making $\hat{\eta}$ a bit smaller, there exists an $\eta>0$ such that, for any $\hat{u}_0\in H^k$ with $\|\hat{u}_0\|_{H^k}\leq K$ and for any $\hat{t}_0\in]T-\hat{\eta},T+\hat{\eta}[$, the solution to the approximate symmetric hyperbolic equation (\ref{eq:symmhypepsilon}) starting in $\hat{u}_0$ at time $\hat{t}_0$ exists on $[\hat{t}_0,\hat{t}_0+\eta[$, and this \emph{independently on $\varepsilon>0$} (use again Step 3).
Taking $\check{\eta}:=\min(\eta,\hat{\eta})>0$, we can look at the initial condition $u(T-\frac{\check{\eta}}{2})$ at time $T-\frac{\check{\eta}}{2}$ and obtain the existence of a family of approximate solutions starting in $u(T-\frac{\check{\eta}}{2})$ at time $T-\frac{\check{\eta}}{2}$ and existing on $[T-\frac{\check{\eta}}{2},T+\frac{\check{\eta}}{2}[$.
Restricting to any compact interval in $[T-\frac{\check{\eta}}{2},T+\frac{\check{\eta}}{2}[$ and applying the preceding results from Steps 4 to 6, we obtain the existence of a solution to the symmetric hyperbolic system starting in $u(T-\frac{\check{\eta}}{2})$ at time $T-\frac{\check{\eta}}{2}$ and existing beyond $T$.
By uniqueness of solutions to symmetric hyperbolic systems, the latter solution coincides with the former on $[T-\frac{\check{\eta}}{2},T[$ and in particular $u$ can be extended beyond $T$, QED.
%The estimate (lower bound) for the time of existence of a local solution in Picard-Lindel\"of can be taken out of mildsolutions1.
\qed\\
\end{proof}

Now we need an additional control on the lifetime of the solution under the additional assumptions of semilinearity (instead of merely quasilinearity) and the one of punctured nonlinearity, i.e., we assume that there is one regular solution (satisfied in our case, as the nonlinearity vanishes at the zero section):

\btheo[Estimate on lifetime]\label{t:lifetime}
Consider a symmetric hyperbolic system of equations on $\mathbb{T}^n$ of the form {\rm(\ref{eq:symmhypsys1storder})} where $A_0,A_j,g\in C^k$ for some $k>\frac{n}{2}+1$.
% \begin{eqnarray}
% \label{eq:symmhypsys1storder}
% A_0(t,x,u)\cdot\frac{\partial u}{\partial t} =  L(t, x, u,\partial)u  + g(t,x, u).
% \end{eqnarray}
\v Assume {\rm (\ref{eq:symmhypsys1storder})} to be {\bf semilinear}, i.e., that $A_0$ and $A_j$ are constant in their last argument $u$, and furthermore assume that there is a sufficiently regular (say, $C^0H^k$) global solution $v$ to {\rm (\ref{eq:symmhypsys1storder})}.\\
Then for every $C,T>0$, there exists an $\varepsilon>0$ such that every $C^1$-solution $u$ to {\rm (\ref{eq:symmhypsys1storder})} with $u(0) =f\in H^k$ and $\|f-v(0)\|_{H^k}\leq\varepsilon$ exists on $[0,T]\times\mathbb{T}^n$ and satisfies $\|u(s)-v(s)\|_{H^k}\leq C$ for all $s\leq T$.
%\item {\bf (stability around zero)} For each interval $I$, there is some $C>0$ such that for all $C^0H^k$-solutions $u$ to {\rm (\ref{basicPDE})} whose $C^0$-norm does not exceed $1$ on $I$ we have $\|u(t) \|_{H^k}\leq e^{Ct} \| f\|_{H^k}$.
\etheo

\v {\bf Remark:} On the one hand, if $g(t,x,0) =0$ for all $(t,x) \in \R\times\mathbb{T}^n$, obviously $0$ is a smooth solution. On the other hand, by defining $\tilde{g} (w) := g(v+w) + (L-A_0\partial_t)v$ for a solution $v$, one can consider the equation $\tilde{P} (w) = 0$ for $\tilde{P} := -A_0\partial_t+L + \tilde{g}$. Obviously, $\tilde{P}(w)=0$ is a symmetric hyperbolic equation, where the nonlinearity $\tilde{g}$ satisfies $\tilde{g} (0) =0$.\\

%into the form $\tilde{L } + \tilde{g}$ for $\tilde{L}$ linear and $\tilde{g} (0) =0$.{\red\bf The problem is that $\tilde{L}$ is no more symmetric hyperbolic as above!}\\

%\footnote{Even if there is a regular non-equilibrium solution $u$ of infinite lifetime one can argue this way but the growth of $u$ enters the estimates, such that a fixed multiple of the logarithm cannot serve any more as a universal estimate, but still for a given lifetime $T$ there is an $\e >0$ such that for all initial values $I$ with $\vert I \vert_{H^k} <\e$ have a lifetime larger than $T$. Same for 'quasilinear' {\bf (?)}.}

\begin{proof} 
%{\red\bf What happens in case the equilibrium solution does not vanish? Because then Moser does not give the estimate below. Adapt the proof.} 
In view of the extension criterion in Theorem \ref{t:exsymmhypsys1storder} and the remark above, we assume that $g(t,x,0)=0$ and estimate $\|u(t)\|_{H^k}^2$ by a function of $t$.
We proceed as in the proof of Claim 3 above and first estimate $\frac{d}{dt}\|J_\varepsilon u(t)\|_{H^k,0}^2$ for any $\varepsilon>0$, where $\|v\|_{H^k,0}^2:=\sum_{|\alpha|\leq k}\left(A_0\cdot \partial^\alpha v,\partial^\alpha v\right)_{L^2}$ for every $v\in H^k$ (both norms $\|\cdot\|_{H^k,0}$ and $\|\cdot\|_{H^k}$ are equivalent on any compact subset of $\R\times\mathbb{T}^n$); then we let $\varepsilon$ tend to $0$ and obtain a differential inequality which, by Bihari's inequality, implies the statement.
Taking into account that all $A_j$, $0\leq j\leq n$, only depend on $(t,x)$ and that $g(t,x,0)=0$, we can mimic the proof of Step 3 and obtain, after letting $\varepsilon\to0$, the estimate
\[\left|\frac{d}{dt}\|u(t)\|_{H^k,0}^2\right|\leq F(\|u(t)\|_{C^0}^2)\cdot\|u(t)\|_{H^k,0}^2\]
for some continuous real-valued function $F$ on $[0,\infty)$.
In particular, up to changing $F$, we obtain 
\[\left|\frac{d}{dt}\|u(t)\|_{H^k,0}^2\right|\leq F(\|u(t)\|_{H^k,0}^2)\cdot\|u(t)\|_{H^k,0}^2.\]
By Bihari's inequality, this proves the statement.
Namely, letting $y(t):=\|u(t)\|_{H^k,0}^2$, we have the inequality $y'\leq yF(y)$ so that, assuming $y>0$ (otherwise $y$ vanishes identically because of Theorem \ref{t:globaluniqnonlinearsymmhyp} below) and setting $z:=\ln(y)$, we obtain
\[\int_{z(0)}^{z(t)}\frac{ds}{F(e^s)}\leq t\]
for every $t\geq0$.
Because $F(e^s)\buil{\longrightarrow}{s\to-\infty}F(0)\geq0$, we have $\int_{z(0)}^{z(t)}\frac{ds}{F(e^s)}\buil{\longrightarrow}{y(0)\searrow0}\infty$ which implies that, for any $T,D>0$, there exists an $\varepsilon>0$ such that, for any $y$ fulfilling $y'\leq yF(y)$ with $y(0)<\varepsilon$, the function $z(t)$ exists on $[0,T]$ and satisfies $z(t)\leq D$.
This concludes the proof.
\findemo
\end{proof}

Finally, we need (global) uniqueness for solutions to symmetric hyperbolic systems.

\btheo\label{t:globaluniqnonlinearsymmhyp}
Consider a $\mathbb{K}^N$-valued first-order symmetric hyperbolic system on $\mathbb{T}^n$ or $\mathbb{R}^n$ as in {\rm Definition \ref{d:symmhypsys1storder}} and assume $A_j,g \in C^1$.
Let $I$ be an open interval with $0\in\,I$.
Let $u_1,u_2\in C^1(I\times\mathbb{T}^n)$ (resp. $u_1,u_2\in C^1(I\times\mathbb{R}^n)$) be any solutions to 
\[A_0(t,x,u_j)\frac{\partial u_j}{\partial t}=L(t,x,u_j,\partial)u_j+g(t,x,u_j)\textrm{ for }t\in I\]
with $u_j(0)=f\in C^0$.
Then $u_1=u_2$.
\etheo 

%Here the space variable $x$ can lie in $\R^n$ or $\mathbb{T}^n$.\\

\begin{proof}
We show that $u_1-u_2$ solves a \emph{linear} symmetric hyperbolic system. 
We write 
\begin{eqnarray*} 
\frac{\partial(u_1-u_2)}{\partial t}&=&A_0^{-1}(t,x,u_1)\cdot L(t,x,u_1,\partial)u_1+A_0^{-1}(t,x,u_1)\cdot g(t,x,u_1)\\
&&-A_0^{-1}(t,x,u_2)\cdot L(t,x,u_2,\partial)u_2-A_0^{-1}(t,x,u_2)\cdot g(t,x,u_2)\\
&=&A_0^{-1}(t,x,u_1)\cdot L(t,x,u_1,\partial)(u_1-u_2)+A_0^{-1}(t,x,u_1)\cdot L(t,x,u_1,\partial)u_2-A_0^{-1}(t,x,u_2)\cdot L(t,x,u_2,\partial)u_2\\
&&+(A_0^{-1}g)(t,x,u_1)-(A_0^{-1}g)(t,x,u_2)\\
&=&A_0^{-1}(t,x,u_1)\cdot L(t,x,u_1,\partial)(u_1-u_2)+(A_0^{-1}(t,x,u_1)-A_0^{-1}(t,x,u_2))\cdot L(t,x,u_1,\partial)u_2\\
&&+A_0^{-1}(t,x,u_2)\cdot(L(t,x,u_1,\partial)-L(t,x,u_2,\partial))u_2+(A_0^{-1}g)(t,x,u_1)-(A_0^{-1}g)(t,x,u_2).
\end{eqnarray*}
Now, because $A_0,A_j,g\in C^1$, we may write 
\begin{eqnarray*} 
A_0^{-1}(t,x,u_1)-A_0^{-1}(t,x,u_2)&=&M(t,x,u_2,u_2)\cdot(u_1-u_2)\\
L(t,x,u_1,\partial)-L(t,x,u_2,\partial)&=&\sum_{j=1}^nB_j(y,x,u_1,u_2)\cdot(u_1-u_2)\frac{\partial}{\partial x_j}\\
(A_0^{-1}g)(t,x,u_1)-(A_0^{-1}g)(t,x,u_2)&=&N(t,x,u_1,u_2)\cdot(u_1-u_2),
\end{eqnarray*} 
therefore 
\begin{eqnarray*} 
\frac{\partial(u_1-u_2)}{\partial t}&=&A_0^{-1}(t,x,u_1)\cdot L(t,x,u_1,\partial)(u_1-u_2)+M(t,x,u_2,u_2)\cdot(u_1-u_2)L(t,x,u_1,\partial)u_2\\
&&+A_0^{-1}(t,x,u_2)\cdot\sum_{j=1}^nB_j(y,x,u_1,u_2)\cdot(u_1-u_2)\frac{\partial u_2}{\partial x_j}+N(t,x,u_1,u_2)\cdot(u_1-u_2),
\end{eqnarray*} 
that is,
\[A_0(t,x,u_1)\cdot\frac{\partial(u_1-u_2)}{\partial t}=L(t,x,u_1,\partial)(u_1-u_2)+B(t,x,u_1,u_2)\cdot(u_1-u_2),\]
where $B$ is of zero order.
Hence $u_1-u_2$ solves a linear symmetric hyperbolic system of first order with vanishing initial condition along the Cauchy hypersurface $\{0\}\times\mathbb{T}^n$ (resp. $\{0\}\times\mathbb{R}^n$) of the globally hyperbolic spacetime $I\times\mathbb{T}^n$ (resp. $I\times\mathbb{R}^n$).
An elementary energy estimate for such systems (see e.g. \cite[Theorem 5.3]{BaerGreenhyperb}) implies that $u_1-u_2=0$ on $I\times\mathbb{T}^n$ (resp. $I\times\mathbb{R}^n$).\hfill\qed\\
\end{proof}

Now we want to transfer the previous local results to the framework of Lorentzian manifolds.
Let $(M^n,g)$ be any globally hyperbolic spacetime and $S\subset M$ be any spacelike Cauchy hypersurface with induced Riemannian metric $g_S$.
Let $E\bui{\longrightarrow}{\pi}M$ be any vector bundle.
A {\bf differential operator $P$ of order $k\in\mathbb{N}$} on $\pi$ is a fibre-bundle-morphism from the $k$th jet bundle $J^k\pi$ of $\pi$ to $\pi$.
It is called {\bf semilinear} if $[\ldots[P,f\cdot],f\cdot,\ldots,f\cdot]=:\sigma_P(df)$ is a vector bundle endomorphism for all scalar functions $f$ on $M$, where $f$ appears $k$ times in the brackets.
Generalizing \cite[Definition 5.1]{BaerGreenhyperb} to the nonlinear case, we define a {\bf semilinear symmetric hyperbolic operator} of first order acting $\pi$ as a semilinear first-order-differential operator $P$ acting on sections of $\pi$ such that, denoting by $\sigma_P: T^*M \rightarrow {\rm End} (E)$ its principal symbol, there is an (definite or indefinite) inner product $\langle\cdot\,,\cdot\rangle$ on $E$ such that for any $\xi\in T^*M$, the endomorphism $\sigma_P(\xi)$ of $E$ is symmetric/Hermitian and positive-definite in case $\xi$ is future-directed causal.
It is easy to see that, locally, $P$ is described exactly by Definition \ref{d:symmhypsys1storder}, where $t$ is a local time-function on $M$.
Theorems \ref{t:exsymmhypsys1storder} and \ref{t:globaluniqnonlinearsymmhyp} imply the following

% Following \cite{BaerGreenhyperb}, we define a {\bf linear symmetric hyperbolic operator} of first order acting on sections of $E\bui{\longrightarrow}{\pi} M$ as a linear first-order-differential operator $P$ acting on sections of $E$ such that, denoting by $\sigma_P$ its principal symbol, there is an (definite or indefinite) inner product $\langle\cdot\,,\cdot\rangle$ on $E$ such that for any $\xi\in T^*M$, the endomorphism $\sigma_P(\xi)$ of $E$ is symmetric/Hermitian and positive-definite in case $\xi$ is future-directed causal.
% It is easy to see that, locally, $P$ is described exactly by Definition \ref{d:symmhypsys1storder}, where $t$ is a local time-function on $M$.
% Theorems \ref{t:exsymmhypsys1storder} and \ref{t:globaluniqnonlinearsymmhyp} imply the following

\bcoro\label{c:existuniqnlinsymmhypglobhyp}
Let $(M^n,g)$ be any globally hyperbolic spacetime and $S\subset M$ be any spacelike Cauchy hypersurface with induced Riemannian metric $g_S$.
Let $E\bui{\longrightarrow}{\pi}M$ be any vector bundle with (definite or indefinite) inner product and $P$ be any semilinear symmetric hyperbolic operator of first order acting on sections of $\pi$.
Let $k\in\mathbb{N}$ with $k>\frac{n-1}{2}+1$.%let $h\colon E\to E$ be any $C^k$ map with $\pi\circ h=\pi$.
Then for any $f\in H^{k,2}(S,g_{S})$, there exists an open neighbourhood $U$ of $S$ in $M$ such that a unique solution $u\in \Gamma_{C^1}(U,E)$ to $Pu=0$ with $u_{|_{S}}=f$ exists.
\ecoro

\begin{proof} Choose for any point $x\in S$ a neighbourhood $B_x$ in $S$ such that the domain of dependence $A_x$ of $B_x$ is contained in a submanifold chart domain for $S$.
Then, via the embedding of $B_x$ into a possibly large torus, we can express the equation $Pu=h(u)$ locally in each $A_x$ as a symmetric hyperbolic system as in Definition \ref{d:symmhypsys1storder}.
Consider for each $x$ a cut-off function which is $1$ on $B_x$ and has support contained in a chart neighbourhood of the torus.
We cut-off the initial data using that function and get the existence of a solution in a small strip around $B_x$.
There is a small neighbourhood of $x$ whose domain of dependence $C_x$ is contained in that strip.
The solutions obtained that way coincide on the intersection of any two such domains.
Patching all such domains $C_x$ together, we obtain a small open neighbourhood of $S$ in $M$ carrying a solution to the original equation.
%{\bf\red Coming soon (Nicolas). }
\hfill\qed
\end{proof}

\bcoro\label{c:lifetime}
Let $(M^n,g)$ be any globally hyperbolic spacetime with compact Cauchy hypersurface $S\subset M$.
Let $k\in\mathbb{N}$ with $k>\frac{n-1}{2}+1$.
Let $E\bui{\longrightarrow}{\pi}M$ be any vector bundle with (definite or indefinite) inner product and $P$ be any $C^k$ semilinear symmetric hyperbolic operator of first order acting on sections of $\pi$ with $P=L+h$, where $L$ is linear and $h$ is of order zero with $h(0)=0$.
Then we have the following {\bf estimate on lifetime} for the solution $u$ of $Pu=0$: for each $T>0$, there is an $\varepsilon>0$ such for all initial values $u_0$ on $S$ with $H^k$-norm smaller than $\varepsilon$, the lifetime for the solution with that initial value is greater than $T$.
\ecoro

\begin{proof}
First observe that for every coordinate patch, a global solution is given by $0$.
Then use finitely many times the estimates given in Theorem \ref{t:lifetime}.\hfill\qed 
\end{proof}

\bigskip

Symmetric hyperbolic operators of second order on $E\bui{\longrightarrow}{\pi}M$ are defined as follows: a differential operator $P$ of second order on $\pi$ is called {\bf symmetric hyperbolic} if there exists a symmetric hyperbolic operator of first order $Q$ -- called {\bf the first prolongation of $P$} -- acting on sections of $\pi\oplus T^*M\!\otimes\!\pi$ such that $Pu=Q(u,\nabla u)$ for every section $u$ of $\pi$. This fits to the restriction to charts --- there, $\nabla u$ is expressed as $\partial u + \Gamma$ where $\Gamma$ is an algebraic (actually, linear) expression in the $u$ variable. Therefore  a representation by $Q$ as above entails an analogous expression in each chart. Furthermore, common textbook knowledge assures that every operator of the form 

$$Pu =   - \partial_t^2 u + \sum_{i,j = 1}^m A_{ij} (t,x) \cdot  \nabla_{ij} u + \sum_{i=1}^m B_i(t,x) \cdot \nabla_i u + c \cdot \partial_t u + d \cdot u     $$

(with $A_{ij}$ symmetric and uniformly positive) can be presented as $Pu = Q(u, \nabla u)$ as above, and the Laplace-d'Alembert equation on a compact subset can be brought into the form $Pu=0$ for $P$ as above.
If $P$ is semilinear, so is $Q$; if $P = P_0 + p$ with $P_0$ linear and $p$ of zeroth order with $p(0) =0$, then $Q= Q_0 + q$ with $Q_0$ linear, $q$ of zeroth order and $q(0)=0$. The local-in-time existence result for second-order symmetric hyperbolic systems is based on Corollary \ref{c:existuniqnlinsymmhypglobhyp}.
It is important to note that, if $P$ has $C^k$ coefficients, then so has $Q$. However, as the new operator $Q$ includes a derivative of $u$, we loose one order of regularity for $u$, but as we do not care much for the weakest possible regularity condition on the initial values anyway, we treat the semilinear operator $Q$ just like a quasilinear operator.
However, notice that there is a folklore theorem mentioned in Taylor's book stating that semilinear symmetric hyperbolic systems of first order have a $C^0$-extension criterion, therefore we could avoid the loss of one derivative of $u$ and obtain sharper statements for the necessary regularity of the initial values.\\

%\end{document}

\begin{small}

\end{small}

\end{document}